\documentclass[11pt,reqno]{amsart}

\usepackage[english]{babel}
\usepackage[utf8]{inputenc}
\usepackage[T1]{fontenc}
\usepackage{lmodern} \normalfont 
\DeclareFontShape{T1}{lmr}{bx}{sc} { <-> ssub * cmr/bx/sc }{}
\usepackage{microtype}	

\usepackage{comment}

\usepackage{amssymb}
\usepackage{amsmath}
\usepackage{amsthm}
\usepackage{mathtools}
\usepackage{mathrsfs}
\usepackage{accents} 
\usepackage{etoolbox}
\usepackage{siunitx}
\usepackage{dsfont} 
\usepackage{graphicx}
\usepackage{color}
\usepackage[dvipsnames]{xcolor}
\usepackage{tikz}
\usetikzlibrary{calc,positioning,shapes}
\usetikzlibrary{patterns,decorations.pathmorphing,decorations.markings}
\usetikzlibrary{spy}
\usepackage{pgfplots}
\pgfplotsset{compat=newest}
\usepackage[margin=10pt,font=small,labelfont=bf,labelsep=endash]{caption}
\usepackage{subcaption}

\usepackage{booktabs}
\usepackage{paralist}
\usepackage{soul}

\numberwithin{equation}{section}

\usepackage{algorithm}
\usepackage{algorithmicx}
\usepackage{algpseudocode}

\usepackage{cite}
\usepackage{multirow} 
\usepackage{enumitem} 
\setlist[enumerate]{label=(\roman*)}

\usepackage[colorlinks,citecolor=black,urlcolor=black]{hyperref}
	\hypersetup{allcolors=black}
\usepackage[nameinlink,noabbrev]{cleveref}
\crefname{subsection}{subsection}{subsections}

\theoremstyle{plain}
\newtheorem{theorem}{Theorem}[section]
\newtheorem{proposition}[theorem]{Proposition}
\newtheorem{lemma}[theorem]{Lemma}
\newtheorem{corollary}[theorem]{Corollary}
\newtheorem{remark}[theorem]{Remark}
\newtheorem{definition}[theorem]{Definition}

\DeclareMathOperator*{\argmin}{arg\,min} 

\definecolor{mycolor1}{rgb}{0.00000,0.44700,0.74100}
\definecolor{mycolor2}{rgb}{0.85000,0.32500,0.09800}
\definecolor{mycolor3}{rgb}{0.92900,0.69400,0.12500}
\definecolor{mycolor4}{rgb}{0.46600,0.67400,0.18800}
\definecolor{mycolor5}{rgb}{0.49400,0.18400,0.55600}

\newcommand{\matlab}{MATLAB\textsuperscript{\textregistered}}

\title{Deterministic Kalman filters for uncertain dynamical systems}

\author{Karl Kunisch${}^{\star,\dagger}$ \and Jesper Schr\"oder${}^\dagger$}
\address{${}^{\star}$ Institute of Mathematics and Scientific Computing, University of Graz, A-8010 Graz, Austria}
\address{${}^{\dagger}$ Johann Radon Institute, Austrian Academy of Sciences, A-4040 Linz, Austria}
\email{karl.kunisch@uni-graz.at}
\email{jesper.schroeder@ricam.oeaw.ac.at}
\date{\today}

\begin{document}
%
\begin{abstract}
	The Kalman(-Bucy) filter is the natural choice for the state reconstruction of disturbed, linear dynamical systems based on flawed and incomplete measurements. Taking a deterministic viewpoint this work investigates possible extensions of the concept to systems with uncertain dynamics and noise covariances. In a theoretical analysis error bounds in terms of the variance of the uncertainties are derived. The article concludes with a numerical implementation of two example systems allowing for a comparison of the estimators. 
\end{abstract}
%

\maketitle


%
\section{Introduction}
%

In this work we aim at reconstructing the state of a linear perturbed dynamical system of the form 
\begin{equation*}
	\begin{aligned}
		\dot{x}(t) &= A_\sigma x(t) + B v(t)
		~~~~~
		t \in (0,T),\\
		x(0) &= x_0 + \eta,
	\end{aligned}
\end{equation*}
where the disturbances in the dynamics and the initial state are given by $v \in L^2(0,T;\mathbb{R}^m)$ and $\eta \in \mathbb{R}^n$, respectively. Additionally the system matrix $A_\sigma$ depends on a parameter $\sigma \in \Sigma$ representing uncertainties of the system dynamics. We assume access to a measurement of the system given as 
\begin{equation*}
	y(t) = C x(t) + \mu(t)
	~~~~~
	t \in (0,T),
\end{equation*}
affected by the disturbance $\mu \in L^2(0,T;\mathbb{R}^r)$.
The disturbances are assumed to be unknown and of deterministic nature. They are associated with parameter dependent positive definite matrices $\Gamma_\sigma \in \mathbb{R}^{n,n}$, $R_\sigma \in \mathbb{R}^{m,m}$, and $Q_\sigma \in \mathbb{R}^{r,r}$ representing their magnitudes. 

Reconstructing a system's state based on incomplete and disturbed measurements has a rich history going back at least to Wiener's seminal work \cite{Wie49}, where disturbances are modeled as random processes. The pioneering works of Kalman and Bucy introduced the highly celebrated Kalman(-Bucy) filter for time discrete \cite{Kal60} and for time continuous \cite{KalB61} linear systems with Gaussian noise. 
Due to its effectiveness it is widely applicable in practice and it inspired numerous variations to address particularly challenging systems such as the extended Kalman filter for nonlinear dynamical systems and the ensemble Kalman filter \cite{Eve94} for systems of very large dimension. 

The assumption of normally distributed noise, however, implies state independent disturbances. Hence exact knowledge of the system dynamics is required and modeling uncertainties and errors introduced by, e.g., material dependent parameters, moving sensors, or model reduction, are not addressed explicitly.
Similarly, the Kalman filter relies on exact knowledge on the distribution of the errors, i.e., their respective covariances which, in practice, are known only approximately. 
In \cite{PetMcf94}, see also the monograph \cite{PetSav99}, the authors consider time continuous, quadratically stable dynamics subject to bounded uncertainties. They design an estimator based on solutions of associated Riccati equations and analyze the asymptotic behavior with respect to time. It is shown that the expected difference of estimator and disturbed state converges to a neighborhood of zero. 
For a discussion of time discrete systems affected by noise with uncertain covariances we refer to \cite{Meh70} and the references therein, see also \cite{ShiJohMu07} for a more recent analysis. 

The remainder of the article is organized as follows. 
\Cref{sec: KF} gives an overview on the Kalman filter in the stochastic as well as in the deterministic setting. We compare the two formulations, point out connections, and set up the framework for our investigation.
In \Cref{sec: KF for uncertain systems} the uncertainties are formulated mathematically and three designs for Kalman filters for uncertain systems are proposed and characterized. In part we shall focus on the estimator which arises as the weighted average of the individual Kalman filters and minimizes the average energy. 
\Cref{sec: errAn} contains the technical proof that the distances of the proposed estimators to the Kalman filter associated with the hidden parameter can be estimated from above in terms of the level of variance in the uncertain matrices.  
The article concludes with a numerical investigation of the designs at the hand of two examples in \Cref{sec: numExp}.

Throughout the article we use the following notation. 
We say that a multivariate random variable $X$ is distributed according to $\mathcal{N}(m,\Sigma)$, where $m \in \mathbb{R}^d$ and $\Sigma \in \mathbb{R}^{d,d}$ is symmetric and positive definite, if $X$ is normally distributed with mean $m$ and covariance matrix $\Sigma$. In short we write $X \sim \mathcal{N}(m,\Sigma)$.
For a symmetric, positive definite matrix $P\in \mathbb{R}^{d,d}$ the associated norm of $x \in \mathbb{R}^d$ is defined via $\Vert x \Vert_M^2 = x^\top M x$. The cardinality of a finite set $\mathcal{M}$ is denoted by $\vert \mathcal{M} \vert$. Finally, we abbreviate some common $L^p$ spaces, Sobolev spaces, and spaces of continuous functions equipped with their default norms as follows. We set $L^2(t_1,t_2;\mathbb{R}^d) = \mathcal{L}_{t_1}^{t_2}$, $H^1(t_1,t_2;\mathbb{R}^d) = \mathcal{H}_{t_1}^{t_2}$, 
$C([t_1;t_2];\mathbb{R}^d) = \mathcal{C}_{t_1}^{t_2} $, and
$C([0;t];\mathbb{R}^{d,d}) = \mathcal{C}_t^{d,d} $. For the spaces of vector valued function the dimension $d$ is suppressed in the notation and becomes clear from context. 

%
\section{The Kalman filter}\label{sec: KF}
%
We commence with a brief presentation of two different formulations of the Kalman filter. Even though our proposed Kalman filter for uncertain systems relies only on the deterministic formulation we include a recap of the stochastic formulation enriching the interpretation of our approach. For ease of presentation we restrict ourselves to time-invariant systems. While the concepts and results presented in this section are well known in the literature we include them for the purpose of a self-contained work. 
%
\subsection{The stochastic formulation}\label{subsec: stochForm}
%
We illustrate the stochastic filter based on the original work \cite{KalB61}. Even though this formulation of the model lacks some mathematical rigor, the formal representation suffices to illustrate the underlying concept. For a more thorough treatment based on It\^o integrals we refer to \cite{Oks98,Xio08}.

Consider the disturbed system with system state $X_t$ and measured output $Y_t$ modeled by
\begin{equation}\label{eq: stochSys}
\begin{aligned}
    \tfrac{\mathrm{d}}{\mathrm{d}t} X_t &= A X_t + B V_t
    ~~~~~ t \in (0,T),\\
    Y_t &= C X_t + W_t,
    ~~~~~ t \in (0,T),
\end{aligned}
\end{equation}
where $A \in \mathbb{R}^{n,n}$ describes the linear dynamics of the system, $C \in \mathbb{R}^{r,n}$ describes the incomplete measurements taken from the state and $B \in \mathbb{R}^{n,m}$ encodes how the noise enters the dynamics. The white noise affecting the dynamics and the measurement is modeled by the multivariate Gaussian random variables $V_t \in \mathbb{R}^m$ and $W_t \in \mathbb{R}^r$, respectively. They have zero mean and covariances given by 
\begin{equation*}
\begin{aligned}
    \text{cov}[V_t,V_s] &= R \, \delta(t-s),\\
    \text{cov}[W_t,W_s] & = Q \, \delta(t-s),
\end{aligned}
\end{equation*}
with given symmetric, positive definite matrices $R \in \mathbb{R}^{m,m}$ and $Q\in \mathbb{R}^{r,r}$. The Gaussian distribution of the state at time $t$ is described by $X_t$ with initial Gaussian distribution $X_0$ with mean $x_0 \in \mathbb{R}^n$ and covariance $\Gamma \in \mathbb{R}^{n,n}$, symmetric, positive definite. Finally, $V_t$, $W_t$, and $X_0$ are assumed to be mutually independent and $Y_t$ encodes the distribution of the output of the system at time $t$.

The objective, as stated in \cite{KalB61}, is to optimally estimate the realized state $x(t)$, of the system based on realized, measured outputs $y(s)$, $0 \leq s \leq t$. Kalman and Bucy find that the optimal estimate $\widehat{x}(t)$ is characterized via
\begin{equation}\label{eq: stochFil}
\begin{aligned}
    \dot{\widehat{x}}(t) &= A \widehat{x}(t) + \Pi(t) C^\top Q^{-1} (y(t) - C\widehat{x}(t))
    ~~~~~t \in (0,T),\\
    \widehat{x}(0) &= x_0,
\end{aligned}
\end{equation}
where $\Pi$ is the unique solution to the differential Riccati equation 
\begin{equation}\label{eq: stochRicc}
\begin{aligned}
    \dot{\Pi}(t) &= A \Pi(t) + \Pi(t) A^\top - \Pi(t) C^\top Q^{-1} C \Pi(t) + B R B^\top
    ~~~~~ t \in (0,T),\\
    \Pi(0) &= \Gamma.
\end{aligned}
\end{equation}
These formulas have since been confirmed by the stochastic community using rigorous mathematical concepts, see e.g. \cite{Oks98,Xio08}. Interpreting \eqref{eq: stochSys} as stochastic differential equations with noise terms defined by Brownian motions and solutions defined via the It\^o integral one characterizes the conditional expectation and covariance. More precisely, given a measured output the state of the system at time $t$ is distributed according to $\mathcal{N}(\widehat{x}(t),\Pi(t))$ where the mean and covariance are characterized via \eqref{eq: stochFil} and \eqref{eq: stochRicc}, respectively. In particular, if one uses the mean $\widehat{x}$ as a state estimator the error has zero mean and a covariance characterized via the Riccati equation \eqref{eq: stochRicc}. 
%
%

We conclude the subsection by introducing some well-known concepts connected to multivariate normal distributions as presented,
\begin{definition}\label{def: concStoch}
    Consider a random variable distributed according to $\mathcal{N}(m,\Sigma)$ with mean $m \in \mathbb{R}^d$ and covariance $\Sigma \in \mathbb{R}^{d,d}$ and assume $\Sigma$ to be regular. We define the following concepts associated with it.
    \begin{enumerate}
        \item The precision matrix \cite[Sec.~5.4]{DeGr70} $P$ is defined as the inverse of the covariance matrix, i.e., $P = \Sigma^{-1}$.
        \item The generalized variance \cite[Sec.~5]{Gir96} is defined as the determinant of the covariance matrix. Further the generalized precision $p$ is defined as the reciprocal of the generalized variance i.e., $p = \frac{1}{\det(\Sigma)} = \det(P)$.
        \item For any fixed vector $\xi \in \mathbb{R}^d$ its \textit{Mahalanobis distance} \cite{Mah36} to the normal distribution is defined as its distance to the mean weighted by the precision, i.e., 
        \begin{equation*}
            \Vert \xi - m \Vert_P = \sqrt{(\xi-m)^\top P (\xi-m)}.
        \end{equation*}
    \end{enumerate}
\end{definition}
The precision matrix and generalized precision extend the precision associated with univariate normal distributions to the multivariate case and the Mahalanobis distance introduces a distance from a point to a point to a given normal distribution. As described above, the Kalman filter can be understood as a time dependent normal distribution. Hence these concepts enable us to quantify the precision of the Kalman filter as well as the distance of deterministic points to the Kalman filter. 

%
\subsection{The deterministic formulation}
%
In this section we consider a deterministic model of a disturbed system brought forward by Mortensen in \cite{Mor68}. In his work he considers possibly nonlinear systems and derives a state estimator based on energy minimization. Unfortunately the obtained formulas for nonlinear systems rely on access to the value function arising from the optimal control problem minimizing the energy. While said value function is formally available as the solution of a Hamilton-Jacobi-Bellman equation its realization is afflicted by the curse of dimensionality. However, Mortensen notes that, when applied to linear dynamics, the formulas reduce to the Kalman filter derived in \cite{KalB61}. A detailed study of this concept for linear systems including a rigorous derivation of the filter equations is presented by Willems in \cite{Wil04}.

The deterministic formulation of the dynamics and the output of the model reads
\begin{equation}\label{eq: detModel}
\begin{alignedat}{2}
    x(t) &= A x(t) + B v(t)
    ~~~~~ &&t \in (0,T),\\
    x(0) &= x_0 + \eta, &&\\
    y(t) &= C x(t) + \mu(t)
    ~~~~~ &&t \in (0,T),
\end{alignedat}
\end{equation}
where the system matrices $A$, $B$, and $C$ are the same as in the stochastic model \eqref{eq: stochSys}. The disturbances in the dynamics, initial value, and output are represented by $v \in L^2(0,T;\mathbb{R}^m)$, $\eta \in \mathbb{R}^n$, and $\mu \in L^2(0,T;\mathbb{R}^r)$, respectively. They are assumed to be deterministic and unknown. For the sake of readability throughout the rest of this work for $0 \leq t_1 < t_2 \leq T $ and $d \in \mathbb{R}^d$ we denote $\mathcal{L}_{t_1}^{t_2} = L^2(t_1,t_2;\mathbb{R}^d)$. Note that the dimension of the image space is surpressed in the notation and is explicitly mentioned if not clear from context.
\begin{remark}\label{rem: L2gap}
    On a formal level system \eqref{eq: detModel} can be interpreted as a version of the stochastic model \eqref{eq: stochSys} in which the random terms have realized but are unknown. As Willems notes, due to a gap in regularity this connection can unfortunately not be established in a rigorous fashion. The assumed $L^2$-regularity of the errors $v$ and $\mu$ is essential for the following analytical treatment. In contrast, realizations of the white noise in the stochastic model (interpreted as the derivative of a Brownian motion) must be expected to be of significantly lower regularity, cf.~\cite{KemScBy24}.
\end{remark}
Willems motivates the construction of the state estimator $\widehat{x}(t)$ at time $t \in (0,T]$ as follows. Given measured output data up until time $t$, i.e., $y(s)$, $0 \leq s \leq t$, one reconstructs the disturbances $\eta$, and $v(s)$, and $\mu(s)$, $0 \leq s \leq t$ as the ones that have minimal (weighted) energy among all possible disturbances that fit the given output in the sense of \eqref{eq: detModel} restricted to $(0,t)$. Mathematically speaking, for a fixed $t \in (0,T]$ the state estimation at time $t$ is defined via 
\begin{equation}\label{eq: energDist}
\begin{alignedat}{2}
    \min_{\eta \in \mathbb{R}^n, v \in \mathcal{L}_0^t, \mu \in \mathcal{L}_0^t}~ 
    &\frac{1}{2} \Vert \eta \Vert_{\Gamma^{-1}}^2 
    + \int_0^t \Vert v(s) \Vert_{R^{-1}}^2 
    + \Vert \mu(s) \Vert_{Q^{-1}}^2 \, \mathrm{d}s,
    ~~~~~
    &&\text{subject to}\\
    x(s) &= A x(s) + B v(s)
    ~~~~~ &&s \in (0,t),\\
    x(0) &= x_0 + \eta, &&\\
    y(s) &= C x(s) + \mu(s)
    ~~~~~ &&s \in (0,t),
\end{alignedat}
\end{equation}
where $\Gamma \in \mathbb{R}^{n,n}$, $R\in \mathbb{R}^{m,m}$, and $Q \in \mathbb{R}^{r,r}$ are symmetric positive definite matrices and in the following are referred to as weighting matrices.
Denoting the resulting state trajectory associated with the minimizing disturbances by $x_t^*$ the estimator at time $t$ is defined as $\widehat{x}(t) = x_t^*(t)$. We stress that for $0 \leq t_1 < t_2 \leq T$ in general it holds 
$\widehat{x}(t_1) 
= x^*_{t_1}(t_1)
\neq x^*_{t_2}(t_1)$.
Therefore, in order to obtain the state estimate $\widehat{x}(t)$ for all $t \in (0,T]$ the minimization must be carried out individually for each time point. Before discussing the derivation of the Kalman filter equations we comment on the weighting matrices $\Gamma$, $R$, and $Q$.
\begin{remark}\label{rem: weighting matrices}
    While the weighting matrices can, in principal, be chosen freely they correspond to the covariance matrices $\Gamma$, $R$, and $Q$ appearing in the stochastic formulation. 
    The intuitive connection is very straight forward if the covariance matrices are diagonal. If, e.g., $Q$ is diagonal the $r$ components of the stochastic measurement error $W_t$ in \eqref{eq: stochSys} are uncorrelated and their respective variances are found on the diagonal of $Q$. Weighting the (distance of the) corresponding deterministic measurement error $\mu$ (to its zero mean) in \eqref{eq: energDist} with $Q^{-1}$ has the following effect. Any component with a small variance (relative to other components of $\mu$ and the components of $\eta$ and $v$) is weighted with the reciprocal, i.e., is strongly penalized in the minimization and the minimizer is promoted to display a small value in said component. This exactly fits the small variance of the stochastic formulation. 
    It is well-known that this concept transfers to the case of non-diagonal covariance matrices $\Gamma$, $R$, and $Q$ allowing for the components of each error to be correlated and it is applied, e.g., in generalized least squares \cite{Ait36}.
    To stress this principle, for the remainder of this work we refer to $\Gamma$, $R$, and $Q$ as covariance rather than weighting matrices.
\end{remark}
In \cite{Wil04} Willems treats the minimization \eqref{eq: energDist} directly. While his work considers the slightly less general setting with fixed $Q = \mathrm{Id}$ and $R = \mathrm{Id}$ a straight forward extension of his arguments leads to analogous results for the setting presented here. One finds that the state estimator $\widehat{x}$ is characterized via the Kalman filter equations \eqref{eq: stochFil}-\eqref{eq: stochRicc}, where now $x_0$ denotes the undisturbed initial condition rather than the mean of the initial distribution. 

For our purposes we follow Mortensen's approach and reformulate \eqref{eq: energDist} into a problem of optimal control. More precisely, we insert the identities for $x(0)$ and $y$ into the cost and additionally fix a final state $\xi \in \mathbb{R}^n$ for the state trajectory to obtain 
\begin{alignat}{2}
    &\min_{x \in \mathcal{H}_0^t, v \in \mathcal{L}_0^t}~ J(x,v;t,\xi) &&= 
    \frac{1}{2} \Vert x(0) - x_0 \Vert_{\Gamma^{-1}}^2 
    + \frac{1}{2} \int_0^t \Vert v(s) \Vert_{R^{-1}}^2 
    + \Vert y - C x(s) \Vert_{Q^{-1}}^2 \, \mathrm{d}s,\label{eq: cost}\\
    &\text{subject to}
    ~~~~~~~~~~\dot{x}(s) &&= A x(s) + B v(s)
    ~~~~~ s \in (0,t), \label{eq: state}\\
    &~~~~~~~~~~~~~~~~~~~~~~~x(t) &&= \xi. \label{eq: FV}
\end{alignat}
\begin{remark}
    We point out that, even though the construction of the estimator is based on the optimal control problem \eqref{eq: cost}-\eqref{eq: FV} and $v$ takes on the role of a control, we do not assume any possibility to influence the system. The control formulation at hand is merely a tool to reconstruct energy minimal disturbances that most likely caused the measured output. 
\end{remark}
The value function associated with the control problem is defined as
\begin{equation}\label{eq: VF}
\begin{aligned}
    \mathcal{V}(t,\xi) &= \inf_{x \in \mathcal{H}_0^t, v \in \mathcal{L}_0^t} J(x,v;t,\xi) ~~~\text{subject to}~ \eqref{eq: state}-\eqref{eq: FV} ~~~ t \in (0,T],\\
    \mathcal{V}(0,\xi) &= \frac{1}{2} \Vert \xi - x_0 \Vert_{\Gamma^{-1}}^2,
\end{aligned}
\end{equation}
and represents the minimal amount of energy required to ensure that at time $t$ the system is in the state $\xi$. It turns out that for all $t \in [0,T]$ the state estimator $\widehat{x}$ defined above can equivalently be defined as
\begin{equation*}
    \widehat{x}(t) = \argmin_{\xi \in \mathbb{R}^n} \mathcal{V}(t,\xi).
\end{equation*}
We now show that the value function can be given explicitly in terms of the solution of the Kalman filter equations which are 
\begin{equation}\label{eq: KF_fixed}
	\begin{aligned}
		\dot{\widehat{x}}(t) 
		&= A \widehat{x}(t)
		+ \Pi(t) C^\top Q^{-1} \left( y(t) - C \widehat{x}(t)\right)
        ~~~~~
		t \in (0,T) ,\\
		\widehat{x}(0) &= x_0,
	\end{aligned}
\end{equation}
and
\begin{equation}\label{eq: Ricc_fixed}
	\begin{aligned}
		\dot{\Pi}(t) &= A \Pi(t) + \Pi(t) A^\top
		- \Pi(t) C^\top Q^{-1} C \Pi(t) 
		+ B R B^\top
		~~~
		~~~t \in (0,T),\\
		\Pi(0) &= \Gamma.
	\end{aligned}
\end{equation}
Even though the following derivations and arguments are based on the deterministic perspective, we still refer to the solution $\Pi$ of \eqref{eq: Ricc_fixed} and its inverse $\Pi^{-1} = P$ as the (error) covariance and the precision matrix associated with the Kalman filter allowing for an interpretation of our results. Multiplying \eqref{eq: Ricc_fixed} by $P(t)$ from the left and from the right we find that the precision matrix is characterized via the differential Riccati equation

\begin{equation}\label{eq: Prec_fixed}
\begin{aligned}
    \dot{P}(t) &= -A^\top P(t) - P(t) A - P(t) BRB^\top P(t) + C^\top Q^{-1} C,
    ~~~
	~~~t \in (0,T),\\
    P(0) &= \Gamma^{-1}.
\end{aligned}
\end{equation}
\begin{lemma}
    Let $\widehat{x}$ and $P$ be given as the solutions of \eqref{eq: KF_fixed} and \eqref{eq: Prec_fixed}, respectively. Then for all $t \in [0,T]$ and $\xi \in \mathbb{R}^n$ it holds
    \begin{equation}\label{eq: VF_fixed}
    	\mathcal{V}(t,\xi) = 
    	\frac{1}{2} (\xi - \widehat{x} (t))^\top P(t) (\xi - \widehat{x} (t))
    	+ \frac{1}{2} \int_0^t \Vert y(s) - C \widehat{x} (s) \Vert_{Q^{-1}}^2 \, \mathrm{d}s.
    \end{equation}
    In particular the value function is quadratic in $\xi$ implying that its Hessian is independent of $\xi$ and its Hessian $\nabla_{\xi\xi}^2 \mathcal{V}(t,\xi)$ is given by the solution $P$ of \eqref{eq: Prec_fixed}.
\end{lemma}
\begin{proof}
    The result is shown in \cite[Sec.~2.3]{BreKu21} for continuous outputs $y$ and weighting matrices $\Gamma = \mathrm{Id}$, $R = \mathrm{Id}$, and $Q = \mathrm{Id}$. A straight forward extension yields the result for the weighted setting. A standard sensitivity analysis, e.g., via an application of the inverse mapping theorem to the first order optimality system ensures that $\mathcal{V}(t,\xi)$ depends continuously on $y \in \mathcal{L}_0^t$. A density argument completes the proof.
\end{proof}

This representation of the value functions allows for an interesting interpretation of the energy. 
\begin{corollary}\label{cor: energyVsDist}
    Let $\widehat{x}$ and $P$ be given as the solution of \eqref{eq: KF_fixed} and \eqref{eq: Prec_fixed}, respectively. Then for all $t \in [0,T]$ and $\xi \in \mathbb{R}^n$ it holds
    \begin{equation*}
        \mathcal{V}(t,\xi) - \mathcal{V}(t,\widehat{x}(t))
        =
        \frac{1}{2} \Vert \xi - \widehat{x}(t) \Vert_{P(t)}^2.
    \end{equation*}
\end{corollary}
%
%
In view of the discussion in \Cref{subsec: stochForm} we find that the difference in energies of a state $\xi$ at time $t$ and the Kalman trajectory $\widehat{x}(t)$ coincides with their squared Mahalanobis distance.
%
\section{Kalman filtering for uncertain systems}\label{sec: KF for uncertain systems}
%
The Kalman filter described in the previous section assumes exact knowledge of both the system dynamics represented by the matrix $A$ and the nature of the occurring noise represented by the weighting/covariance matrices $\Gamma$, $R$, and $Q$. In practical applications, however, this information may not be available in exact form as the dynamics of the system may depend on various uncertain parameters \cite{ShiJohMu07}. Further, as described, e.g. in \cite{Meh70}, the covariance matrices associated with the distribution of the noise and the initial state may be subject to uncertainties. 
%
\subsection{Modeling uncertain parameters}\label{subsec: ModUnc}
%
We model the uncertainties of the system matrix $A$ as follows. 
Let $s_A \in \mathbb{N}$ and consider a matrix valued function
\begin{equation*}
    a \colon \mathbb{R}^{s_A} \rightarrow \mathbb{R}^{n,n}.
\end{equation*}
Further let 
$N_A \in \mathbb{N}$
and 
$\Sigma_A \subset \mathbb{R}^{s_A}$
be of cardinality 
$\vert \Sigma_A \vert = N_A$. Consider the discrete probability space $(\Sigma_A, \mathcal{P}(\Sigma_A),P_A)$
, where 
$\mathcal{P_A}(\Sigma_A)$
denotes the power set of $\Sigma_A$ and for $S \subset \Sigma_A$ we define the probability function 
$P_A(S) = \frac{1}{N_A} \vert S \vert$.
Then the linear dynamics of the system are characterized by the matrix valued random variable 
\begin{equation*}
    \sigma \mapsto a_\sigma = a(\sigma),~~~ \Sigma_A \to \mathbb{R}^{n,n},
\end{equation*}
and $a_\sigma$ denotes the system matrix that realizes with the parameter $\sigma \in \Sigma_A$.
We analogously define the uncertain weighting matrices $\gamma_\sigma$, $r_\sigma$, and $q_\sigma$ as random variables on $\Sigma_\Gamma$, $\Sigma_R$, and $\Sigma_Q$ with respective cardinalities and dimensions $N_\Gamma$, $s_\Gamma$, $N_R$, $s_R$, and $N_Q$, $s_Q$. The associated probabilities are defined and denoted analogously to $A$.
We denote $N = N_A \, N_\Gamma \, N_R \, N_Q$ and $s = s_A \, s_\Gamma \, s_R \, s_Q$. Then
\begin{equation*}
    \Sigma = \Sigma_A \times \Sigma_\Gamma \times \Sigma_R \times \Sigma_Q \subset \mathbb{R}^s
\end{equation*}
has cardinality $N$ and
we define the product space $(\Sigma, \mathcal{F},P)$, via
\begin{equation*}
\begin{aligned}
    \mathcal{F} &= \mathcal{P}(\Sigma_A) \otimes \mathcal{P}(\Sigma_\Gamma) \otimes \mathcal{P}(\Sigma_R) \otimes \mathcal{P}(\Sigma_Q) = \mathcal{P}(\Sigma),\\
    P(S_A \otimes S_\Gamma \otimes S_R \otimes S_Q) &= P_A(S_A) \, P_\Gamma(S_\Gamma) \, P_R(S_R)\, P_Q(S_Q) = \frac{1}{\vert S_A \otimes S_\Gamma \otimes S_R \otimes S_Q \vert}.
\end{aligned}
\end{equation*}
Finally the uncertain system and weighting matrices we work with are defined as random variables on $\Sigma$. Denote $\sigma = (\sigma_A,\sigma_\Gamma,\sigma_R,\sigma_Q) \in \Sigma$ and define
\begin{equation*}
\begin{alignedat}{2}
    \Sigma &\to \mathbb{R}^{n,n},
    &&\sigma  \mapsto A_\sigma =  a_{\sigma_A},\\
    \Sigma &\to \mathbb{R}^{n,n},
    &&\sigma  \mapsto \Gamma_\sigma =  \gamma_{\sigma_\Gamma},\\
    \Sigma &\to \mathbb{R}^{m,m},
    ~~~&&\sigma  \mapsto R_\sigma = r_{\sigma_R},\\
    \Sigma &\to \mathbb{R}^{r,r},
    &&\sigma  \mapsto Q_\sigma = q_{\sigma_Q}.
\end{alignedat}
\end{equation*}
Note that they are, by construction, mutually independent. 
For convenience we denote the $N$ elements of $\Sigma$ by $\sigma_k$, $k = 1,\dots,N$ and dependencies on $\sigma_k$ may be indicated by indexing $k$, e.g., $A_{\sigma_k} = A_k$.

In the following we present three different estimators for the reconstruction of the state of the uncertain system affected by noise with uncertain covariances. 
%
\subsection{The Kalman filter for expected system and covariance matrices}\label{subs: KF for expected system}
%
We begin with the most straight forward and computationally affordable approach. Namely, we form the expectation of the four uncertain matrices and realize the associated Kalman filter. 
More precisely, define 
\begin{equation*}
\begin{alignedat}{2}
    A_\mathbb{E} &= \mathbb{E}[A_\sigma] = \frac{1}{N} \sum_{k=1}^N A_k,
    ~~~~~
    \Gamma_\mathbb{E} &&= \mathbb{E}[\Gamma_\sigma] = \frac{1}{N} \sum_{k=1}^N \Gamma_k,\\
    R_\mathbb{E} &= \mathbb{E}[R_\sigma] = \frac{1}{N} \sum_{k=1}^N R_k,
    ~~~~~
    Q_\mathbb{E} &&= \mathbb{E}[Q_\sigma] = \frac{1}{N} \sum_{k=1}^N Q_k.
\end{alignedat}
\end{equation*}
The associated state estimator is constructed as the solution to
\begin{equation}\label{eq: expected obs}
\begin{aligned}
    \dot{\widehat{x}}_\mathbb{E}(t) 
    &= A_\mathbb{E} \widehat{x}_\mathbb{E}(t)
    + \Pi_\mathbb{E}(t) C^\top Q_\mathbb{E}^{-1} \left( y(t) - C \widehat{x}_\mathbb{E}(t)\right)
    ~~~~~
    t \in (0,T) ,\\
    \widehat{x}(0) &= x_0,
\end{aligned}
\end{equation}
where $\Pi_\mathbb{E}$ is given via the Riccati equation
\begin{equation}\label{eq: expected Ricc}
\begin{aligned}
    \dot{\Pi}_\mathbb{E}(t) &= A_\mathbb{E} \Pi_\mathbb{E}(t) + \Pi_\mathbb{E}(t) A_\mathbb{E}^\top
    - \Pi_\mathbb{E}(t) C^\top Q_\mathbb{E}^{-1} C \Pi_\mathbb{E}(t) 
    + B R_\mathbb{E} B^\top
    ~~~
    ~~~t \in (0,T),\\
    \Pi(0) &= \Gamma_\mathbb{E}.
\end{aligned}
\end{equation}
%

%
\subsection{The expected Kalman filter}\label{subs: expKal}
%
For the second approach we realize the Kalman filter associated with each possible parameter $\sigma \in \Sigma$ allowing for the construction of an estimator $\widehat{x}_\sigma$ depending on the uncertainty. Formulating the uncertain trajectory as a random variable its associated expectation yields the state estimator. 

To that end for every element $\sigma_k \in \Sigma$ solve for the associated Kalman filter $\widehat{x}_k = \widehat{x}_{\sigma_k}$ via 
\begin{equation}\label{eq: MemKF}
\begin{aligned}
    \dot{\widehat{x}}_k(t) 
    &= A_k \widehat{x}_k(t)
    + \Pi_k(t) C^\top Q_k^{-1} \left( y(t) - C \widehat{x}_k(t)\right)
    ~~~~~
    t \in (0,T) ,\\
    \widehat{x}_k(0) &= x_0,
\end{aligned}
\end{equation}
where $\Pi_k$ is given via the Riccati equation
\begin{equation}\label{eq: MemRicc}
\begin{aligned}
    \dot{\Pi}_k(t) &= A_k \Pi_k(t) + \Pi_k(t) A_k^\top
    - \Pi_k(t) C^\top Q_k^{-1} C \Pi_k(t) 
    + B R_k B^\top
    ~~~
    ~~~t \in (0,T),\\
    \Pi_k(0) &= \Gamma_k.
\end{aligned}
\end{equation}
For later reference we denote the associated precision matrix by $P_k(t) \coloneqq \Pi_k^{-1}(t)$. It is given as the unique solution of 
\begin{equation}\label{eq: MemPrec}
\begin{aligned}
    \dot{P}_k(t) &= -A_k^\top P_k(t) - P_k(t) A_k - P_k(t) BR_k B^\top P_k(t) + C^\top Q_k^{-1} C,
    ~~~
	~~~t \in (0,T),\\
    P_k(0) &= \Gamma_k^{-1}.
\end{aligned}
\end{equation}
Now $\sigma_k \mapsto \widehat{x}_{\sigma_k}$, $ \Sigma \to \{ \widehat{x}_{\sigma_k} \colon \sigma_k \in \Sigma \} \subset \mathcal{H}_0^T$ is considered a random variable where the image space is equipped with its power set as a $\sigma$-algebra.
The state estimator then is defined as the expectation, i.e.,
\begin{equation}\label{eq: def_varnothing}
    \widehat{x}_\varnothing \coloneqq \mathbb{E}[\widehat{x}_{\sigma}] = \frac{1}{N} \sum_{k=1}^N \widehat{x}_k.
\end{equation}
We point out that, by construction, for any $t \in [0,T]$ the estimator $\widehat{x}_{\varnothing}(t)$ minimizes the expected squared Euclidean distance, i.e., it is the unique solution of 
\begin{equation*}
    \min_{\xi \in \mathbb{R}^n} \mathbb{E} \left[ \Vert \xi - \widehat{x}_\sigma(t) \Vert^2 \right]
    = \min_{\xi \in \mathbb{R}^n} \frac{1}{N} \sum_{k=1}^N \Vert \xi - \widehat{x}_k(t) \Vert^2.   
\end{equation*}
Finally, note that the realization of this option is computationally more demanding than the one in the previous subsection. Indeed, it requires solving $N$ Riccati equations and ODEs. 

%
\subsection{Minimizing the expected energy}\label{subs: MinEner}
%
The third and final estimator is defined as the minimizer of the expected energy. 
Throughout the rest of this work we use the notation established in the previous subsection, i.e., $\widehat{x}_k$, $\Pi_k$, and $P_k$ are the Kalman filter trajectory, the error covariance, and the precision matrix associated with the matrices $A_k$, $\Gamma_k$, $R_k$, and $Q_k$. Additionally, throughout the remainder of this work we denote the associated value function according to \eqref{eq: VF} by $\mathcal{V}_k = \mathcal{V}_{\sigma_k}$. For $t \in [0,T]$ we obtain
\begin{equation*}
    \widehat{x}_k(t) = \argmin_{\xi \in \mathbb{R}^n} \mathcal{V}_k(t,\xi),
\end{equation*}
and
\begin{equation}\label{eq: VF_k}
    \mathcal{V}_k(t,\xi) = 
    \frac{1}{2} (\xi - \widehat{x}_k (t))^\top P_k(t) (\xi - \widehat{x}_k (t))
    + \frac{1}{2} \int_0^t \Vert y(s) - C \widehat{x}_k (s) \Vert_{Q_k^{-1}}^2 \, \mathrm{d}s,
\end{equation}
with $\widehat{x}_k$ and $P_k$ given via \eqref{eq: MemKF} and \eqref{eq: MemPrec}, respectively.
Note that analogously to $\widehat{x}_k= \widehat{x}_{\sigma_k}$ the value function $\mathcal{V}_k = \mathcal{V}_{\sigma_k}$ can be considered as a random variable. The reader is further reminded that intuitively speaking $\mathcal{V}_k(t,\xi)$ represents the minimum amount of energy required to ensure that the system associated with $\sigma_k$ at time $t$ is in the state $\xi$. 

For the uncertain system we now define the expected energy as follows.

\begin{definition}\label{def: AverEn} 
	The expected energy $E \colon [0,T] \times \mathbb{R}^n \to \mathbb{R}$ is defined via
	\begin{equation}\label{eq: AverEn}
    \begin{aligned}
		E(t,\xi) 
        &= \mathbb{E}\left[ \mathcal{V}_\sigma(t,\xi) \right]
        = \frac{1}{N} \sum_{k=1}^N \mathcal{V}_k(t,\xi)\\
		&= \frac{1}{2N} \sum_{k=1}^N \left(
        (\xi - \widehat{x}_k (t))^\top P_k(t) (\xi - \widehat{x}_k (t))
        + \int_0^t \Vert y(s) - C \widehat{x}_k (s) \Vert_{Q_k^{-1}}^2 \, \mathrm{d}s
        \right).
    \end{aligned}
	\end{equation}
\end{definition}
Intuitively, $E(t,\xi)$ represents the amount of energy that is required in expectation to ensure that at time $t$ the uncertain system is in the state $\xi$. This motivates the formal definition of the third state estimator discussed in this work.
\begin{definition}\label{def: UQKF}
	For $t \in [0,T]$ we define the state estimator as the state minimizing the expected energy via
	\begin{equation*}
		\widehat{x}_\mathrm{E}(t) = \argmin_{\xi \in \mathbb{R}^n} E(t,\xi).
	\end{equation*}
\end{definition}
In the following we show that for each $t \in [0,T] $ the function $E(t,\cdot)$ admits a unique minimizer ensuring that $\widehat{x}_\mathrm{E}$ is in fact well-defined. We further derive a characterization in terms of the individual Kalman trajectories $\widehat{x}_k = \widehat{x}_{\sigma_k}$ allowing for a numerical realization with cost comparable to the realization of $\widehat{x}_\varnothing$.

The representation via the individual value functions shows that $E$ is quadratic in $\xi$ allowing for straight forward differentiation of $E$ with respect to $\xi$.
\begin{lemma}
	For any fixed $t \in [0,T]$ the mapping $\xi \mapsto E(t,\xi)$ is of class $C^\infty$. Its gradient and Hessian are given as
	\begin{equation*}
		\nabla E(t,\xi) = \frac{1}{N} \sum_{k=1}^N P_k(t) (\xi - \widehat{x}_k(t)),
		~~~~~
		\text{and}
		~~~~~
		\nabla^2 E(t,\xi) = \frac{1}{N}\sum_{k=1}^N P_k(t),
	\end{equation*}
	respectively.
\end{lemma}
\begin{proof}
	According to \Cref{def: AverEn} the average energy is a sum of functions quadratic in $\xi$ and differentiability follows. The formulas are obtained by differentiation with respect to $\xi$.  
\end{proof}
Existence, uniqueness, and the characterization of a minimizer follow directly.
\begin{proposition}\label{prop: char_en_min}
	The state estimator defined in \Cref{def: UQKF} is well-defined and given in terms of the Kalman filter trajectories and covariance matrices associated with the individual Kalman filters. For $t \in [0,T]$ it holds
	\begin{equation}\label{eq: en_min}
		\widehat{x}_\mathrm{E}(t) = \mathcal{P}^{-1}(t)
							\sum_{k=1}^N P_k(t) \, \widehat{x}_k(t), 
	\end{equation}
    where $\mathcal{P}(t) = \sum_{k=1}^N P_k(t) $.
\end{proposition}
\begin{proof}
	We fix $t \in [0,T]$ and begin by showing uniqueness.
	Due to the regularity of $E$ established in the previous lemma any minimizer $\xi^*$ must satisfy $\nabla E(t,\xi^*) = 0$. Rearranging the terms in the formula for the gradient yields
	\begin{equation*}
		\sum_{k=1}^N P_k(t) \, \xi^* = \sum_{k=1}^N P_k(t) \, \widehat{x}_k(t).
	\end{equation*}
	For all $k = 1,\dots,N$ we have that $P_k$ is given as the solution to a differential Riccati equation with positive definite initial value, hence $P_k(t)$ is positive definite for all $t \in [0,T]$, cf.~\cite[Prop.~1.1]{DieEi94}. The positive definiteness then carries over to the to the sum denoted by $\mathcal{P}$. Hence $\mathcal{P}$ is invertible and $\widehat{x}_\mathrm{E}(t)$ is the only state satisfying the necessary condition to be a minimizer. Together with the formula for the Hessian these considerations additionally ensure positive definiteness of the Hessian and it is shown that $\widehat{x}(t)$ is in fact the unique minimizer.
\end{proof}
\begin{remark}
    According to this formula $\widehat{x}_\mathrm{E}$ can be interpreted as a generalized weighted mean of the individual Kalman filter trajectories. Indeed, it is given as the sum of each Kalman trajectory multiplied from the left by its associated precision matrix. The sum is then multiplied by the inverse of the sum of those precision matrices.
    We note that the expression for $\widehat x_{\mathrm{E}}$ as a weighted average of the family $\widehat x_k$ resembles the statistical method of inverse-variance weighting  were a family of  random variables is aggregated  with the goal to minimize the variance of the weighted average. 
\end{remark}
Interestingly $\widehat{x}_\mathrm{E}$ not only minimizes the expected energy, but also the expected squared Mahalanobis distance.
\begin{corollary}\label{cor: energyMin}
    Let $t \in [0,T]$ be arbitrary. Then $\widehat{x}_\mathrm{E}(t)$ is the unique solution of the minimization problem 
    \begin{equation*}
        \min_{\xi \in \mathbb{R}^n} \mathbb{E} \left[ \Vert \xi - \widehat{x}_\sigma(t) \Vert_{P_\sigma(t)}^2 \right]
        = \min_{\xi \in \mathbb{R}^n} \frac{1}{N} \sum_{k=1}^N \Vert \xi - \widehat{x}_k(t) \Vert_{P_k(t)}^2.
    \end{equation*}
\end{corollary}
\begin{proof}
    The assertion is a direct consequence of \Cref{cor: energyVsDist} and \Cref{def: AverEn}.
\end{proof} 
Before we turn to the error analysis we mention a fourth estimator. 
\begin{remark}
    Another possible design for the reconstruction of the state of the uncertain system is obtained by solving the Kalman filter equation for the mean \eqref{eq: expected obs} where for the gain one replaces the solution $\Pi_\mathbb{E}$ of the associated Riccati equation \eqref{eq: expected Ricc} by the expected gain, i.e., 
    $\frac{1}{N} \sum_{k=1}^N \Pi_k$. It can be understood as an intermediate option combining the designs of \Cref{subs: KF for expected system} and \Cref{subs: expKal}. 
\end{remark}
%

%
\section{Error analysis}\label{sec: errAn}
%
In this section we discuss the error that is introduced by the uncertainty in system and noise covariances. Our analysis is based on the following premise. If one had access to the true, hidden parameter $\bar{\sigma} \in \Sigma$ the obvious choice for the state estimator would be the associated Kalman filter $\widehat{x}_{\bar{\sigma}}$ with error covariance
$\Pi_{\bar{\sigma}}$ given via
\begin{equation}\label{eq: trueKal}
\begin{aligned}
    \dot{\widehat{x}}_{\bar{\sigma}}(t) 
    &= A_{\bar{\sigma}} \widehat{x}_{\bar{\sigma}}(t)
    + \Pi_{\bar{\sigma}}(t) C^\top Q_{\bar{\sigma}}^{-1} \left( y(t) - C \widehat{x}_{\bar{\sigma}}(t)\right)
    ~~~~~
    t \in (0,T) ,\\
    \widehat{x}_{\bar{\sigma}}(0) &= x_0,
\end{aligned}
\end{equation}
and
\begin{equation}\label{eq: Ricc_true}
\begin{aligned}
    \dot{\Pi}_{\bar{\sigma}}(t) &= A_{\bar{\sigma}} \Pi_{\bar{\sigma}}(t) + \Pi_{\bar{\sigma}}(t) A_{\bar{\sigma}}^\top
    - \Pi_{\bar{\sigma}}(t) C^\top Q_{\bar{\sigma}}^{-1} C \Pi_{\bar{\sigma}}(t) 
    + B R_{\bar{\sigma}} B^\top
    ~~~
    ~~~t \in (0,T),\\
    \Pi_{\bar{\sigma}}(0) &= \Gamma_{\bar{\sigma}}.
\end{aligned}
\end{equation}
Hence the quality of the three estimators presented in \Cref{sec: KF for uncertain systems} is measured in terms of their respective distances to $\widehat{x}_{\bar{\sigma}}$.
More precisely, for a given estimator $\widehat{x}$ and time $t \in [0,T]$ we analyze the distance
\begin{equation}\label{eq: errSig}
    \Vert \widehat{x}(t) - \widehat{x}_{\bar{\sigma}}(t) \Vert_{\Pi_{\bar{\sigma}}^{-1}(t)}
    = \Vert \widehat{x}(t) - \widehat{x}_{\bar{\sigma}}(t) \Vert_{P_{\bar{\sigma}}(t)}.
\end{equation}
The justification for the particular choice of weighted norm lies in the stochastic interpretation. In this case the estimator of the state at time $t$ associated with the true parameter $\bar{\sigma}$ is given by a normal distribution with mean $\widehat{x}_{\bar{\sigma}}(t)$ and covariance $\Pi_{\bar{\sigma}}(t)$. The weighting then is motivated by the reasoning illustrated in \Cref{rem: weighting matrices}. 
We further note that in the spirit of \Cref{cor: energyVsDist} it holds that
\begin{equation*}
    \frac{1}{2}\Vert \widehat{x}(t) - \widehat{x}_{\bar{\sigma}}(t) \Vert_{P_{\bar{\sigma}}(t)}^2
    =
    \mathcal{V}_{\bar{\sigma}}(t,\widehat{x}(t))
    - \mathcal{V}_{\bar{\sigma}}(t,\widehat{x}_{\bar{\sigma}}(t)).
\end{equation*}
Hence by estimating the weighted difference we simultaneously estimate the surplus of energy when compared to the minimal energy. 

To facilitate a concise presentation of the results we introduce the following notation. For a given parameter $\sigma$ the associated four-tuple consisting of system matrix and covariance matrices is denoted by $\mathcal{S}_\sigma = (A_{{\sigma}},\Gamma_{{\sigma}}, R_{{\sigma}},Q_{{\sigma}}) \in \mathbb{R}^{n,n} \times \mathbb{R}^{n,n} \times \mathbb{R}^{m,m} \times \mathbb{R}^{r,r}$. For $p \in \mathbb{N}$ we introduce the norm
$
\Vert \mathcal{S}_{\sigma} \Vert_p^p
=
\Vert A_\sigma  \Vert_{\mathbb{R}^{n,n}}^p
+ \Vert \Gamma_\sigma \Vert_{\mathbb{R}^{n,n}}^p
+ \Vert R_\sigma \Vert_{\mathbb{R}^{m,m}}^p
+ \Vert Q_\sigma \Vert_{\mathbb{R}^{r,r}}^p
$ 
and note that sums of systems are understood componentwise.

%
\subsection{Technical preparations}\label{subs: TechPrep}
%
Consider $\widehat{x}_*$ and $\Pi_*$ given as the unique solutions of
\begin{equation}\label{eq: Obs_star}
\begin{aligned}
    \dot{\widehat{x}}_*(t) 
    &= A_* \widehat{x}_*(t)
    + \Pi_*(t) C^\top Q_*^{-1} \left( y(t) - C \widehat{x}_*(t)\right)
    ~~~~~
    t \in (0,T) ,\\
    \widehat{x}_*(0) &= x_0,
\end{aligned}
\end{equation}
and
\begin{equation}\label{eq: Ricc_star}
\begin{aligned}
    \dot{\Pi}_*(t) &= A_* \Pi_*(t) + \Pi_*(t) A_*^\top
    - \Pi_*(t) C^\top Q_*^{-1} C \Pi_*(t) 
    + B R_* B^\top
    ~~~
    ~~~t \in (0,T),\\
    \Pi_*(0) &= \Gamma_*.
\end{aligned}
\end{equation}
These systems will be used with  $(\widehat{x}_*,\Pi_*) = (\widehat{x}_\mathbb{E},\Pi_\mathbb{E})$ and $(\widehat{x}_*,\Pi_*) = (\widehat{x}_k,\Pi_k)$  
to obtain error estimates for the three state estimators presented in \Cref{sec: KF for uncertain systems}.
In order to derive an estimate for $\Vert \Pi_* - \Pi_{\bar{\sigma}} \Vert_{\mathcal{C}_T^{n,n}}$ we first link the covariance matrices to the optimal control problems dual to the state estimation problem. 
For this purpose we define $\Sigma_*(t) = \Pi_*(T-t)$ for $t \in [0,T]$.  By differentiating $t \mapsto \Pi_*(T-t)$ and utilizing the equation for $\dot{\Pi}_*$ we find that $\Sigma_*$ is the unique solution to 
\begin{equation*}
\begin{aligned}
    \dot{\Sigma}_*(t) &= - A_* \Sigma_*(t) - \Sigma_* A_*^\top
                        + \Sigma_*(t) C^\top Q_*^{-1} C \Sigma_*(t) - B R_* B^\top
                        ~~~~~ t \in (0,T),\\
    \Sigma_*(T) &= \Gamma_*. 
\end{aligned}
\end{equation*}
This is the Riccati equation associated to the optimal control problem
\begin{align}
    &\min_{x \in \mathcal{H}_t^T, u \in \mathcal{L}_t^T}~ \mathcal{J}_*(x,u;t,x_0) = 
    \frac{1}{2} \int_t^T \Vert u(s) \Vert_{Q_*}^2 
    + \Vert B^\top x(s) \Vert_{R_*}^2 \, \mathrm{d}s
    +\frac{1}{2} \Vert x(T) \Vert_{\Gamma_*}^2,\label{eq: OCP_star_cost} \\
    \begin{split}
    &\text{subject to}
    ~~~~~~~~~~\dot{x}(s) = A_*^\top x(s) + C^\top u(s)
    ~~~~~ s \in (t,T), \label{eq: OCP_star}\\
    &~~~~~~~~~~~~~~~~~~~~~~~x(t) = x_0,
    \end{split}
\end{align}
for $t \in [0,T)$.
In the following lemma we recall that the associated value function 
\begin{equation*}
    \nu_*(t,x_0) = \inf_{x,u \text{ s.t.} \eqref{eq: OCP_star}  } \mathcal{J}_*(x,u;t,x_0),
\end{equation*}
can be represented via the solution of the associated Riccati equation, cf.~\cite[Thm.~37]{Son98}.
\begin{lemma}
    Let $\Sigma_*$  be as introduced above. Then
    for all $t \in [0,T)$ and $x \in \mathbb{R}^n$ it holds
    \begin{equation*}
        \nu_*(t,x) =\frac{1}{2} \left( x, \Sigma_*(t) x \right)
    \end{equation*}
\end{lemma}
We also define $\Sigma_{\bar{\sigma}}$ and $J_{\bar{\sigma}}$ associated with the state equation  $\dot{x}(s) = A^\top_{\bar \sigma} x(s) + C^\top u(s)$ and obtain an analogous representation of the related value function $\nu_{\bar{\sigma}}$ in terms of $\Sigma_{\bar{\sigma}} = 
\Pi_{\bar{\sigma}}(T-\cdot)$. 
Utilizing these characterizations we can represent the difference of covariance matrices via the difference of value functions.
\begin{lemma}\label{lem: DiffCov_eq_DiffVal}
    Let $\Pi_{\bar{\sigma}}$ and $\Pi_*$ be the solutions of \eqref{eq: Ricc_true} and \eqref{eq: Ricc_star}, respectively. Then for all $t \in (0,T]$ it holds
    \begin{equation*}
        \frac{1}{2} \Vert \Pi_*(t) - \Pi_{\bar{\sigma}}(t) \Vert
        = \sup_{\Vert x \Vert = 1} \left\vert \nu_*(T-t,x) - \nu_{\bar{\sigma}}(T-t,x) \right\vert,
    \end{equation*}
    where $\nu_*$ and $\nu_{\bar{\sigma}}$ are the value functions defined above. 
\end{lemma}
\begin{proof}
    Since both $\Sigma_*$ and $\Sigma_{\bar{\sigma}}$ are symmetric, the same holds for their difference. For all $t \in [0,T]$ it follows 
    \begin{equation*}
    \begin{aligned}
        &\frac{1}{2} \Vert \Pi_*(t) - \Pi_{\bar{\sigma}}(t) \Vert 
        = \frac{1}{2} \Vert \Sigma_*(T-t) - \Sigma_{\bar{\sigma}}(T-t) \Vert 
        = \frac{1}{2} \sup_{\Vert x \Vert = 1} \left\vert \left(x, (\Sigma_*(T-t) - \Sigma_{\bar{\sigma}}(T-t))x \right) \right\vert\\
        &=  \sup_{\Vert x \Vert = 1} \left\vert \frac{1}{2} \left(x, \Sigma_*(T-t) x \right)
            - \frac{1}{2} \left(x, \Sigma_{\bar{\sigma}}(T-t) x \right) \right\vert
        =  \sup_{\Vert x \Vert = 1} \left\vert \nu_*(T-t,x) - \nu_{\bar{\sigma}}(T-t,x) \right\vert,
    \end{aligned}
    \end{equation*}
    and the proof is complete.
\end{proof}
To obtain the desired estimate for the difference of covariances we now estimate the difference of value functions.
First we summarize some estimates for the solution of the optimal control problem associated with $\nu_*$.
\begin{lemma}\label{lem: est*}
    Let $(\bar{x}_*,\bar{u}_*)$ be the unique solution to \eqref{eq: OCP_star_cost}-\eqref{eq: OCP_star} with associated adjoint state $p^*$. There exists $c_*>0$ independent of $t$ such that
    \begin{equation*}
    \begin{aligned}
        \Vert \bar{x}_* \Vert_{\mathcal{H}_t^T}^2
        + \Vert \bar{x}_* \Vert_{\mathcal{C}_t^T}^2
        + \Vert \bar{u}_* \Vert_{\mathcal{L}_t^T}^2
        + \Vert p_* \Vert_{\mathcal{H}_t^T}^2
        + \Vert p_* \Vert_{\mathcal{C}_t^T}^2
        \leq c_* \Vert x_0 \Vert^2.
    \end{aligned}
    \end{equation*}
\end{lemma}
\begin{proof} 
    Since this results can be proven using standard techniques from optimal control, we only sketch the proof. 
    First note that
    \begin{equation}\label{eq: est*}
        \Vert \bar{u}_* \Vert_{\mathcal{L}_t^T}^2
        \leq 2 \mathcal{J}_*(\bar{x}_*,\bar{u}_*;t,x_0)
        \leq 2 \mathcal{J}_*(x^0,0;t,x_0),
    \end{equation}
    where $x^0$ is the unique solution on $(t,T)$ to
    \begin{equation*}
        \dot{x}^0 = A_*^\top x^0,  
        ~~~~~~\dot{x}^0(t) = x_0.
    \end{equation*}
    Using standard arguments including an application of Gronwall's Lemma one can show the existence of $c > 0$ independent of $t$ such that
    $\Vert x^0 \Vert_{\mathcal{H}_t^T}^2 \leq c \Vert x_0 \Vert^2$. Together with \eqref{eq: est*} and the continuous embedding $\mathcal{H}_t^T \hookrightarrow \mathcal{C}_t^T$ with constant $c^\prime$ we get
    \begin{equation*}
    \begin{aligned}
        \Vert \bar{u}_* \Vert_{\mathcal{L}_t^T}^2
        \leq 
        \int_t^T \Vert B x^0 \Vert_{R^*}^2 \, \mathrm{d}s
        + \Vert x^0 (T) \Vert^2
        \leq
        c^2 \Vert B \Vert^2 \, \Vert {R^*}^{\tfrac{1}{2}} \Vert^2 \Vert x_0 \Vert^2
        + 
        c^2 {c^\prime}^2 \Vert x_0 \Vert^2.
    \end{aligned}
    \end{equation*}
    The estimate transfers to $\bar{x}_*$ and $p_*$ using the state and adjoint equations and arguments involving Gronwall's Lemma.
\end{proof}
We now employ the strategy presented in \cite[Lem.~3.2]{GutKuRo24} to show the desired estimate for the value functions.
\begin{lemma}\label{lem: estJ_sigma}
    Let $t \in [0,T)$ and $x_0 \in \mathbb{R}^n$ with $\Vert x_0 \Vert \leq 1$. Then there exists a constant $c_\mathcal{J} > 0 $ independent of $t$ such that 
    \begin{equation*}
        \mathcal{J}_{\bar{\sigma}}
        ( \bar{x}_* - \bar{x}_{\bar{\sigma}},
        \bar{u}_* - \bar{u}_{\bar{\sigma}};
        t, x_0)
        \leq
        c_\mathcal{J}
        \Vert \mathcal{S}_* - \mathcal{S}_{\bar{\sigma}} \Vert_2^2
    \end{equation*}
    where $(\bar{x}_*,\bar{u}_*)$ and $(\bar{x}_{\bar{\sigma}},\bar{u}_{\bar{\sigma}})$ are the optimal pairs associated with $\mathcal{J}_*$ and $\mathcal{J}_{\bar{\sigma}}$, respectively.
\end{lemma}
\begin{proof}
    We have $\nu_*(t,x_0) = \mathcal{J}_*(\bar{x}_*,\bar{u}_*;t,x_0)$, where $(\bar{x}_*,\bar{u}_*)$ is the unique minimizer of \eqref{eq: OCP_star_cost}-\eqref{eq: OCP_star} which is characterized via the optimality system
    \begin{equation*}
    \begin{alignedat}{2}
        \dot{\bar{x}}_*(s) &= A_*^\top \bar{x}_*(s) + C^\top \bar{u}_*(s),~~~~~~~~~~
        &&\bar{x}_*(t) = x_0,\\
        \dot{p}_*(s) &= - A_* p_*(s)- BR_*B^\top \bar{x}_*(s),
        &&p_*(T) = \Gamma_* \bar{x}_*(T),\\
        \bar{u}_*(s) &= - Q_*^{-1} C p_*(s), 
    \end{alignedat}
    \end{equation*}
    where the time dependent equations hold for $s \in (t,T)$. 
    Analogously we have $\nu_{\bar{\sigma}}(t,x_0) = \mathcal{J}_{\bar{\sigma}}(\bar{x}_{\bar{\sigma}},\bar{u}_{\bar{\sigma}};t,x_0)$, where $(\bar{x}_{\bar{\sigma}},\bar{u}_{\bar{\sigma}})$ is the corresponding unique minimizer characterized via the optimality system
    \begin{equation*}
    \begin{alignedat}{2}
        \dot{\bar{x}}_{\bar{\sigma}}(s) &= A_{\bar{\sigma}}^\top \bar{x}_{\bar{\sigma}}(s) + C^\top \bar{u}_{\bar{\sigma}}(s),~~~~~~~~~~
        &&\bar{x}_{\bar{\sigma}}(t) = x_0,\\
        \dot{p}_{\bar{\sigma}}(s) &= - A_{\bar{\sigma}} p_{\bar{\sigma}}(s)- BR_{\bar{\sigma}}B^\top \bar{x}_{\bar{\sigma}}(s),
        &&p_{\bar{\sigma}}(T) = \Gamma_{\bar{\sigma}} \bar{x}_{\bar{\sigma}}(T),\\
        \bar{u}_{\bar{\sigma}}(s) &= - Q_{\bar{\sigma}}^{-1} C p_{\bar{\sigma}}(s).
    \end{alignedat}
    \end{equation*}
    We define
    \begin{equation*}
    \begin{alignedat}{3}
        x_{\delta} &\coloneqq \bar{x}_* - \bar{x}_{\bar{\sigma}},
        ~~~~~
        u_{\delta} &&\coloneqq \bar{u}_* - \bar{u}_{\bar{\sigma}},
        ~~~~~~~~~~
        p_{\delta} &&\coloneqq \bar{p}_* - \bar{p}_{\bar{\sigma}},\\
        \Gamma_{\delta} & \coloneqq \Gamma_* - \Gamma_{\bar{\sigma}},
        ~~~~~
        R_{\delta} &&\coloneqq R_* - R_{\bar{\sigma}},
        ~~~~~
        Q_{\delta^{-1}} &&\coloneqq Q_*^{-1} - Q_{\bar{\sigma}}^{-1},
        ~~~~~
        A_{\delta} = A_* - A_{\bar{\sigma}},
    \end{alignedat}
    \end{equation*}
    and obtain
    \begin{equation*}
    \begin{alignedat}{2}
        \dot{x}_{\delta} 
        &= A_{\bar{\sigma}}^\top x_{\delta}
            + A_{\delta}^\top \bar{x}_*
            + C^\top u_{\delta},
        ~~~~~
        ~~~~~
        ~~~~~
        ~~~~~
        ~~~~~
        ~~~~
        x_{\delta} (t) &&= 0,\\
        \dot{p}_{\delta}
        &= - A_{\bar{\sigma}} p_{\delta} 
            - A_{\delta} p_* 
            - B R_{\bar{\sigma}} B^\top x_{\delta}
            - B R_{\delta} B^\top \bar{x}_*,
        ~~~~~
        p_{\delta} (T) 
        &&= \Gamma_{\bar{\sigma}} x_{\delta} (T)
            + \Gamma_{\delta} \bar{x}_*(T),\\
        u_{\delta} 
        &= - Q_{\bar{\sigma}}^{-1} C p_{\delta} 
            - Q_{\delta^{-1}} C p_*.
    \end{alignedat}
    \end{equation*}
    Now a combination of standard calculus, Gronwall's Lemma, the norm bound on $x_0$, and \Cref{lem: est*} yield existence of a constant $c_\delta>0$ depending on $T$, $c_*$, $\Vert A_{\bar{\sigma}} \Vert$, $\Vert B \Vert$, $\Vert C\Vert$, $\Vert R_{\bar{\sigma}} \Vert$, and $\Vert \Gamma_{\bar{\sigma}} \Vert$ but independent of $t$ such that
    \begin{equation}\label{eq: est_x_p_delta}
    \begin{aligned}
        \Vert x_{\delta} \Vert_{\mathcal{C}_t^T}^2
        +
    	\Vert x_{\delta} \Vert_{\mathcal{H}_t^T}^2
    	&\leq c_\delta \left( \Vert A_\delta \Vert^2 
    	+ 
    	\Vert u_\delta \Vert_{\mathcal{L}_t^T}^2 \right),\\
        \Vert p_\delta \Vert_{\mathcal{C}_t^T}^2
        +
    	\Vert p_\delta \Vert_{\mathcal{H}_t^T}^2	
    	&\leq c_\delta \left(
    	\Vert u_\delta \Vert_{\mathcal{L}_t^T}^2
    	+ 
    	\Vert A_\delta \Vert^2 
    	+ 
    	\Vert R_\delta \Vert^2 
    	+
    	\Vert \Gamma_\delta \Vert^2 
    	\right).
    \end{aligned}
    \end{equation}

    Further, by testing the equation for $\dot{p}_{\delta}$ with $x_{\delta}$ we find 
    \begin{equation}\label{eq: dual1}
        -\int_t^T (\dot{p}_{\delta}, x_{\delta}) \, \mathrm{d}s
        =
        \int_t^T
        (A_{\bar{\sigma}} p_{\delta}, x_{\delta} )
        + 
        (A_{\delta} p_*, x_{\delta} )
        +
        (B R_{\bar{\sigma}} B^\top x_{\delta}, x_{\delta} )
        +
        (B R_{\delta} B^\top \bar{x}_*, x_{\delta})
        \, \mathrm{d}s.
    \end{equation}
    On the other hand, with $x_{\delta} (0) = 0$ and the identity for $p_{\delta} (T)$ we find testing the equation for $\dot x_\delta$ with $p_\delta$
    \begin{equation}\label{eq: dual2}
    \begin{aligned}
        &-\int_t^T (\dot{p}_{\delta}, x_{\delta}) \, \mathrm{d}s
        = \int_t^T (\dot{x}_{\delta}, p_{\delta}) \, \mathrm{d}s
            - (p_{\delta}(T),x_{\delta}(T))\\
        &= \int_t^T 
        (A_{\bar{\sigma}}^\top x_{\delta}, p_{\delta})
        + (A_{\delta}^\top \bar{x}_*, p_{\delta} )
        + ( C^\top u_{\delta}, p_{\delta} )
        \, \mathrm{d}s
        - ( \Gamma_{\bar{\sigma}} x_{\delta}(T), x_{\delta}(T) )
        - ( \Gamma_{\delta} \bar{x}_*, x_{\delta}(T) ).
    \end{aligned}
    \end{equation}
    By subtracting \eqref{eq: dual1} from \eqref{eq: dual2}, rearranging terms and utilizing the identity for $u_{\delta}$ we arrive at
    \begin{equation*}
    \begin{aligned}
        \delta \mathcal{J}_{\bar{\sigma}}
        &\coloneqq
        \int_t^T 
        \Vert u_{\delta} \Vert_{Q_{\bar{\sigma}}}^2
        + \Vert B^\top x_{\delta} \Vert_{R_{\bar{\sigma}}}^2
        \, \mathrm{d}s
        + \Vert x_{\delta}(T) \Vert_{\Gamma_{\bar{\sigma}}}^2\\
        &=
        - \int_t^T  
        (A_{\delta} p_* , x_{\delta})
        + (B R_{\delta} B^\top \bar{x}_*, x_{\delta})
        - ( A_{\delta}^\top \bar{x}_*, p_{\delta} )
        + ( Q_{\bar{\sigma}} u_{\delta}, Q_{\delta^{-1}} p_* )
        \, \mathrm{d}s\\
        &- ( \Gamma_{\delta} \bar{x}_*(T), x_{\delta}(T) ).
    \end{aligned}
    \end{equation*}
    Utilizing \Cref{lem: est*}, \eqref{eq: est_x_p_delta}, and Young's inequality it follows that for any $\epsilon>0$ it holds
    \begin{equation*}
    \begin{aligned}
    	\delta \mathcal{J}_{\bar{\sigma}}
    	&\leq
    	c_* \int_t^T
    		\Vert A_\delta \Vert \Vert x_\delta \Vert
    		+
    		\Vert B \Vert^2 \Vert R_\delta \Vert \Vert x_\delta \Vert
    		+
    		\Vert A_\delta \Vert  \Vert p_\delta \Vert
    		+
    		  \Vert Q_{\bar{\sigma}} \Vert \Vert Q_{\delta^{-1}} \Vert \Vert u_\delta \Vert 
     	\, \mathrm{d}s\\
     	&+
        c_* \Vert \Gamma_\delta \Vert \Vert x_\delta \Vert_{\mathcal{C}_t^T}\\
     	&\leq
     	c_* \int_t^T
     	\frac{1}{\epsilon} \Vert A_\delta \Vert^2
     	+
     	\epsilon \Vert x_\delta \Vert^2
     	+
     	\frac{\Vert B \Vert^4}{2\epsilon} \Vert R_\delta \Vert^2
     	+
     	\frac{\epsilon}{2} \Vert p_\delta \Vert^2 
     	+
     	\frac{\Vert Q_{\bar{\sigma}} \Vert^2}{2 \epsilon} \Vert Q_{\delta^{-1}} \Vert^2
     	+ 
     	\frac{c_* \epsilon}{2} \Vert u_\delta \Vert^2
     	\, \mathrm{d}s\\
     	&+ 
     	\frac{c_*}{2\epsilon} \Vert \Gamma_\delta \Vert^2
     	+   
     	\frac{\epsilon}{2} \Vert x_\delta \Vert_{\mathcal{C}_t^T}^2\\
     	&\leq
     	\frac{c_* T}{\epsilon} \Vert A_\delta \Vert^2
     	+
     	\frac{c_* T \Vert B \Vert^4}{2\epsilon} \Vert R_\delta \Vert^2
     	+
     	\frac{c_* T \Vert Q_{\bar{\sigma}} \Vert^2}{2\epsilon} \Vert Q_{\delta^{-1}} \Vert^2
     	+
     	\frac{c_*}{2\epsilon} \Vert \Gamma_\delta \Vert^2
        + 
     	\frac{c_* \epsilon}{2} \Vert u_\delta \Vert_{\mathcal{L}_t^T}^2\\
     	&+ 
     	(\epsilon c_* c_\delta + \frac{\epsilon c_\delta}{2}) (\Vert A_\delta \Vert^2 + \Vert u_\delta \Vert_{\mathcal{L}_t^T}^2)
     	+ \frac{\epsilon c_* c_\delta}{2} \left(
     	\Vert u_\delta \Vert_{\mathcal{L}_t^T}^2
     	+ 
     	\Vert A_\delta \Vert^2 
     	+ 
     	\Vert R_\delta \Vert^2 
     	+
     	\Vert \Gamma_\delta \Vert^2 
     	\right).
	\end{aligned}  
    \end{equation*}
    Since 
    \begin{equation*}
        \Vert u_\delta \Vert_{\mathcal{L}_t^T}^2 \leq \Vert Q_{\bar{\sigma}}^{-\tfrac{1}{2}} \Vert^2 \int_t^T \Vert u_\delta \Vert_{Q_{\bar{\sigma}}}^2 \, \mathrm{d}s,
    \end{equation*}
    for a sufficiently small $\epsilon > 0$ all terms on the right hand side depending on $u_\delta$ can be absorbed in the left hand side, and with appropriate constants $c_1, c_2 > 0$ it follows that 
    \begin{equation*}
        c_1 \delta \mathcal{J}_{\bar{\sigma}}
        \leq 
        c_2 \left( 
        \Vert A_\delta \Vert + \Vert \Gamma_\delta \Vert + \Vert R_\delta \Vert + \Vert Q_{\delta^{-1}} \Vert
        \right).
    \end{equation*}
    Finally, note that the mapping $ \mathcal{F} (A) = A^{-1}$ is continuously differentiable on
    the open set of positive definite matrices with derivative $D\mathcal{F}(A)B = - A^{-1}B A^{-1}$. Hence an application of Taylor's theorem yields
    \begin{equation*}
        \Vert Q_{\delta^{-1}} \Vert
        =
        \Vert Q_*^{-1} - Q_{\bar{\sigma}}^{-1} \Vert
        \leq 
        \Vert Q_{\bar{\sigma}}^{-1} \Vert^2 \Vert Q_* - Q_{\bar{\sigma}} \Vert
    \end{equation*}
    and the assertion follows.
\end{proof}
With this estimate at hand we are in a position to estimate the difference of value functions. 
\begin{lemma}\label{lem: est_diff_cov}
    Let $\Pi_{\bar{\sigma}}$ and $\Pi_*$ be the solutions of \eqref{eq: Ricc_true} and \eqref{eq: Ricc_star}, respectively. Then there exists a constant $c_\Pi >0$ such that the following inequality holds
    \begin{equation*}
    \begin{aligned}
        \Vert \Pi_* - \Pi_{\bar{\sigma}} \Vert_{\mathcal{C}_t^{n,n}}
        \leq
        c_\Pi 
        \Vert \mathcal{S}_* - \mathcal{S}_{\bar{\sigma}} \Vert_1.
    \end{aligned}
    \end{equation*}
\end{lemma}
\begin{proof}
    Due to \Cref{lem: DiffCov_eq_DiffVal} it suffices to estimate the difference of the value functions. For this purpose let $(\bar{x}_*,\bar{u}_*)$ and $(\bar{x}_{\bar{\sigma}},\bar{u}_{\bar{\sigma}})$ be as in \Cref{lem: estJ_sigma}, let $x_0 \in \mathbb{R}^n$ with $\Vert x_0 \Vert = 1$, and let $t \in [0,T)$ be fixed. It follows
    \begin{equation}\label{eq: estValF}
    \begin{aligned}
        &\left\vert \nu_*(t,x_0) - \nu_{\bar{\sigma}}(t,x_0) \right\vert
        =
        \vert \mathcal{J}_*(\bar{x}_*,\bar{u}_*;t,x_0) 
        -
        \mathcal{J}_{\bar{\sigma}}(\bar{x}_{\bar{\sigma}},\bar{u}_{\bar{\sigma}};t,x_0) \vert\\
        &\leq
        \frac{1}{2} \int_t^T \vert \Vert \bar{u}_* \Vert_{Q_*}^2 - \Vert \bar{u}_{\bar{\sigma}} \Vert_{Q_{\bar{\sigma}}}^2 \vert
        +
        \vert \Vert B^\top \bar{x}_* \Vert_{R_*}^2 - \Vert B^\top \bar{x}_{\bar{\sigma}} \Vert_{R_{\bar{\sigma}}}^2 \vert
        \, \mathrm{d} s
        +
        \frac{1}{2} \vert \Vert \bar{x}_* (T) \Vert_{\Gamma_*}^2 - \Vert \bar{x}_{\bar{\sigma}} (T) \Vert_{\Gamma_{\bar{\sigma}}}^2 \vert.
    \end{aligned}
    \end{equation}
    We estimate the terms on the right hand side individually. For the term involving the controls consider
    \begin{equation*}
    \begin{aligned}
        \int_t^T \vert
        \Vert \bar{u}_* \Vert_{Q_*}^2 
        - \Vert \bar{u}_{\bar{\sigma}} \Vert_{Q_{\bar{\sigma}}}^2 
        \vert \, \mathrm{d}s
        &\leq 
        \int_t^T 
        \vert
            (\bar{u}_*, Q_* \bar{u}_*) - (\bar{u}_*, Q_{\bar{\sigma}} \bar{u}_*)
        \vert
        +
        \vert
            (\bar{u}_*, Q_{\bar{\sigma}} \bar{u}_*) - (\bar{u}_{\bar{\sigma}}, Q_{\bar{\sigma}} \bar{u}_{\bar{\sigma}})
        \vert
        \, \mathrm{d}s\\
        \leq&
        \Vert Q_* - Q_{\bar{\sigma}} \Vert \int_t^T \Vert \bar{u}_* \Vert^2 \, \mathrm{d}s
        +
        \Vert Q_{\bar{\sigma}}^{\frac{1}{2}} \Vert
        \int_t^T
            \Vert \bar{u}_* + \bar{u}_{\bar{\sigma}} \Vert
            \, \Vert \bar{u}_* - \bar{u}_{\bar{\sigma}} \Vert_{Q_{\bar{\sigma}}}
        \, \mathrm{d}s.
    \end{aligned}
    \end{equation*}
    With \Cref{lem: est*} and an analogous result for the solution of the optimal control problem formulated in terms of $\mathcal{J}_{\bar{\sigma}}$ with some constant $c_{\bar{\sigma}}>0$ we obtain
    \begin{equation*}
        \int_t^T \vert
        \Vert \bar{u}_* \Vert_{Q_*}^2 
        - \Vert \bar{u}_{\bar{\sigma}} \Vert_{Q_{\bar{\sigma}}}^2 
        \vert \, \mathrm{d}s
        \leq
        c_* \Vert Q_* - Q_{\bar{\sigma}} \Vert
        +
        (\sqrt{c_*} +  \sqrt{c_{\bar{\sigma}}}) \sqrt{T} \, \Vert Q_{\bar{\sigma}}^{\frac{1}{2}} \Vert \, \sqrt{ \int_t^T \Vert \bar{u}_* - \bar{u}_{\bar{\sigma}} \Vert_{Q_{\bar{\sigma}}}^2 \, \mathrm{d}s }.
    \end{equation*}
    Similarly we find
    \begin{equation*}
    \begin{aligned}
        &\int_t^T \vert
        \Vert B^\top \bar{x}_* \Vert_{R_*}^2 
        - \Vert B^\top \bar{x}_{\bar{\sigma}} \Vert_{R_{\bar{\sigma}}}^2 
        \vert \, \mathrm{d}s\\
        &\leq
        c_* \Vert B \Vert^2 \Vert R_* - R_{\bar{\sigma}} \Vert 
        +
        (\sqrt{c_*} +  \sqrt{c_{\bar{\sigma}}}) \sqrt{T} \, \Vert B \Vert \, \Vert R_{\bar{\sigma}}^{\frac{1}{2}} \Vert
        \sqrt{\int_t^T \Vert B^\top ( \bar{x}_* - \bar{x}_{\bar{\sigma}}) \Vert_{R_{\bar{\sigma}}}^2 \, \mathrm{d}s}
    \end{aligned}
    \end{equation*}
    and
    \begin{equation*}
        \vert \Vert \bar{x}_* (T) \Vert_{\Gamma_*}^2 - \Vert \bar{x}_{\bar{\sigma}} (T) \Vert_{\Gamma_{\bar{\sigma}}}^2 \vert
        \leq
        c_* \Vert \Gamma_* - \Gamma_{\bar{\sigma}} \Vert
        +
        \Vert \Gamma_{\bar{\sigma}}^{\tfrac{1}{2}} \Vert (\sqrt{c_*} + \sqrt{c_{\bar{\sigma}}})
        \Vert \bar{x}_*(T) - \bar{x}_{\bar{\sigma}} (T) \Vert_{\Gamma_{\bar{\sigma}}}.
    \end{equation*}
    Hence with \eqref{eq: estValF} and $c_1 = c_* \max \left( \Vert B \Vert^2,1 \right)$ and an appropriate $c_2 >0$ it follows
    \begin{equation*}
    \begin{aligned}
        \left\vert \nu_*(t,x_0) - \nu_{\bar{\sigma}}(t,x_0) \right\vert
        &\leq 
        c_1 \left(
        \Vert \Gamma_* - \Gamma_{\bar{\sigma}} \Vert
        + 
        \Vert R_* - R_{\bar{\sigma}} \Vert
        +
        \Vert Q_* - Q_{\bar{\sigma}} \Vert
        \right)\\
        &+ c_2 \sqrt{ \mathcal{J}_{\bar{\sigma}}(\bar{x}_* - \bar{x}_{\bar{\sigma}},\bar{u}_* - \bar{u}_{\bar{\sigma}};t,x_0) }.
    \end{aligned}
    \end{equation*}
    Finally with \Cref{lem: DiffCov_eq_DiffVal} and \Cref{lem: estJ_sigma} we obtain
    \begin{equation*}
    \begin{aligned}
        \Vert \Pi_*(T-t) - \Pi_{\bar{\sigma}}(T-t) \Vert
        &\leq c_1 \left(
        \Vert \Gamma_* - \Gamma_{\bar{\sigma}} \Vert
        + 
        \Vert R_* - R_{\bar{\sigma}} \Vert
        +
        \Vert Q_* - Q_{\bar{\sigma}} \Vert
        \right)\\
        &+
        c_2 \sqrt{c_{\mathcal{J}}} 
        \left( 
        \Vert A_* - A_{\bar{\sigma}} \Vert
        + \Vert \Gamma_* - \Gamma_{\bar{\sigma}} \Vert
        + \Vert R_* - R_{\bar{\sigma}} \Vert
        + \Vert Q_* - Q_{\bar{\sigma}} \Vert
        \right),
    \end{aligned}
    \end{equation*}
    implying the estimate for $t \in (0,T]$. Due to the initial condition of the Riccati equations the estimate holds in $t = 0$ as well and the assertion is shown. 
\end{proof}
Combining the estimate for the difference of error covariances with another application of Gronwall's Lemma we estimate the difference of the state reconstruction trajectories. 
\begin{lemma}\label{lem: est_Kal}
    Let $\widehat{x}_*$, $\Pi_*$, $\widehat{x}_{\bar{\sigma}}$, and $\Pi_{\bar{\sigma}}$ be the solutions of \eqref{eq: Obs_star}, \eqref{eq: Ricc_star}, \eqref{eq: trueKal}, and \eqref{eq: Ricc_true} respectively. Then there exists a constant $c_{x_*} > 0$ independent of $t\in[0,T]$ such that
    \begin{equation*}
    \begin{aligned}
        \Vert \widehat{x}_*(t) - \widehat{x}_{\bar{\sigma}}(t) \Vert
        \leq
        c_{x_*} 
        \Vert \mathcal{S}_* - \mathcal{S}_{\bar{\sigma}} \Vert_1.
    \end{aligned}
    \end{equation*} 
\end{lemma}
\begin{proof}
    For $s \in [0,T]$ denote $e(s) = \widehat{x}_*(s) - \widehat{x}_{\bar{\sigma}}(s)$. Subtracting \eqref{eq: trueKal} from \eqref{eq: Obs_star} and rearranging terms shows that $e$ is characterized as the unique solution to
    \begin{equation*}
    \begin{aligned}
        \dot{e}(s) 
        &= \left( A_{\bar{\sigma}} - \Pi_{\bar{\sigma}} (s) C^\top Q_{\bar{\sigma}}^{-1} C \right) e(s)
        + \left( \Pi_*(s) - \Pi_{\bar{\sigma}}(s) \right) C^\top Q_{\bar{\sigma}}^{-1} \left(  y(s) - C \widehat{x}_*(s) \right)\\
        &+ \left( A_* - A_{\bar{\sigma}} \right) \widehat{x}_*(s)
        + \Pi_*(s) C^\top \left( Q_*^{-1} - Q_{\bar{\sigma}}^{-1} \right) (y(s) - C \widehat{x}_*(t)),
                        ~~~~~ s \in (0,T),\\
        e(0) &= 0.
    \end{aligned}
    \end{equation*}
    Testing the first equation with $e(s)$ and integrating over $(0,t)$ for some fixed $t \in(0,T]$ yields
    \begin{equation*}
    \begin{aligned}
        \frac{1}{2} \Vert e(t) \Vert^2
        &\leq \int_0^t \left( \Vert A_{\bar{\sigma}} \Vert + \Vert C \Vert^2 \Vert Q_{\bar{\sigma}}^{-1} \Vert \Vert \Pi_{\bar{\sigma}}(s) \Vert \right) \Vert e (s) \Vert^2
        + \Vert A_* - A_{\bar{\sigma}} \Vert \Vert e(s) \Vert\\
        &+ \Vert \Pi_*(s) - \Pi_{\bar{\sigma}}(s) \Vert \Vert C \Vert \Vert Q_{\bar{\sigma}}^{-1} \Vert 
            \Vert y(s) - C \widehat{x}_*(s) \Vert \Vert e(s) \Vert\\
        &+ \Vert \Pi_*(s) \Vert  \Vert C \Vert \Vert Q_*^{-1} - Q_{\bar{\sigma}}^{-1} \Vert \Vert y(s) - C \widehat{x}_*(s) \Vert \Vert e(s) \Vert 
        \, \mathrm{d}s.
    \end{aligned}
    \end{equation*}
    With Young's inequality we find
    \begin{equation*}
    \begin{aligned}
        \Vert e(t) \Vert^2
        &\leq 
        \int_0^t \left(2 \Vert A_{\bar{\sigma}} \Vert 
        + 2 \Vert C \Vert^2 \Vert Q_{\bar{\sigma}}^{-1} \Vert \Vert \Pi_{\bar{\sigma}}(s) \Vert
        + 3 \right) \Vert e (s) \Vert^2
        \, \mathrm{d}s
        + T \Vert A_* - A_{\bar{\sigma}} \Vert^2\\
        &+ \Vert C \Vert^2
        \Vert Q_{\bar{\sigma}}^{-1} \Vert^2 
        \Vert y - C \widehat{x}_* \Vert_{\mathcal{L}_0^T}^2 
        \Vert \Pi_* - \Pi_{\bar{\sigma}} \Vert_{\mathcal{C}_T^{n,n}}^2\\
        &+ \Vert C \Vert^2 \Vert
        \Pi_* \Vert_{\mathcal{C}_T^{n,n}}^2
        \Vert y - C \widehat{x}_* \Vert_{\mathcal{L}_0^T}^2
        \Vert Q_*^{-1} - Q_{\bar{\sigma}}^{-1} \Vert^2, 
    \end{aligned}
    \end{equation*}
    and finally Gronwall's inequality yields
    \begin{equation*}
        \Vert e(t) \Vert
        \leq
        c \left(
        \Vert A_* - A_{\bar{\sigma}} \Vert 
        + \Vert Q_*^{-1} - Q_{\bar{\sigma}}^{-1} \Vert
        + \Vert \Pi_* - \Pi_{\bar{\sigma}} \Vert_{\mathcal{C}_T^{n,n}}
        \right),
    \end{equation*} 
    for some appropriate constant $c>0$.
    Again we estimate $\Vert Q_*^{-1} - Q_{\bar{\sigma}}^{-1} \Vert$ in terms of $\Vert Q_* - Q_{\bar{\sigma}} \Vert$ as we did in the proof of \Cref{lem: estJ_sigma}.
    The assertion now follows with \Cref{lem: est_diff_cov} and a constant $c_{x_*}>0$
    depending on $T$, 
    $\Vert C \Vert$, $\Vert A_{\bar{\sigma}} \Vert$, $ \Vert Q_{\bar{\sigma}}^{-1} \Vert $, 
    $\Vert y - C \widehat{x}_* \Vert_{\mathcal{L}_0^T} $,
    $\Vert \Pi_* \Vert_{\mathcal{C}_T^{n,n}}$,
    and $\Vert \Pi_{\bar{\sigma}} \Vert_{\mathcal{C}_T^{n,n}}$.
\end{proof}
This concludes the technical preparations and we turn to the error estimation for the three state estimators discussed in \Cref{sec: KF for uncertain systems}.
%
\subsection{Worst case error estimation}
%
We begin by estimating the difference of the estimator $\widehat{x}_{\mathbb{E}}$ introduced in \Cref{subs: KF for expected system}. 
\begin{proposition}\label{prop: est_exp}
    There exists a constant $c_\mathbb{E}>0$ independent of $N$ such that for all $t \in [0,T]$ it holds
    \begin{equation*}
    \begin{aligned}
        \Vert \widehat{x}_{\mathbb{E}}(t) - \widehat{x}_{\bar{\sigma}}(t) \Vert_{P_{\bar{\sigma}(t)}}
        &\leq
        c_\mathbb{E}~
        \mathbb{E} \left[
        \Vert \mathcal{S}_\sigma - \mathcal{S}_{\bar{\sigma}} \Vert_1
        \right].
    \end{aligned}
    \end{equation*}
\end{proposition}
\begin{proof}
    We apply the results of \Cref{subs: TechPrep} with 
    \begin{equation*}
    \begin{aligned}
        A_* &= A_{\mathbb{E}} = \frac{1}{N} \sum_{k=1}^N A_k,
        ~~~~~~
        \Gamma_* = \Gamma_{\mathbb{E}} = \frac{1}{N} \sum_{k=1}^N \Gamma_k,\\
        R_* &= R_{\mathbb{E}} = \frac{1}{N} \sum_{k=1}^N R_k,
        ~~~~~~
        Q_* = Q_{\mathbb{E}} = \frac{1}{N} \sum_{k=1}^N Q_k,
    \end{aligned}
    \end{equation*}
    and \Cref{lem: est_Kal} yields existence of $c_{x_*} > 0$ independent of $N$ such that for any $t \in [0,T]$ it holds
    \begin{equation*}
    \begin{aligned}
        \Vert \widehat{x}_{\mathbb{E}}(t) - \widehat{x}_{\bar{\sigma}}(t) \Vert_{P_{\bar{\sigma}(t)}}
        &\leq 
        \Vert P_{\bar{\sigma}}^{\tfrac{1}{2}}(t) \Vert
        ~
        \Vert \widehat{x}_{\mathbb{E}}(t) - \widehat{x}_{\bar{\sigma}}(t) \Vert\\
        &\leq
        c_{x_*} \,
        \Vert P_{\bar{\sigma}}^{\tfrac{1}{2}}(t) \Vert
        \left(
        \Vert A_{\mathbb{E}} - A_{\bar{\sigma}} \Vert
        +\Vert \Gamma_{\mathbb{E}} - \Gamma_{\bar{\sigma}} \Vert
        +\Vert R_{\mathbb{E}} - R_{\bar{\sigma}} \Vert
        +\Vert Q_{\mathbb{E}} - Q_{\bar{\sigma}} \Vert
        \right).
    \end{aligned}
    \end{equation*}
    We note that $P_{\bar{\sigma}}$ is continuous in $[0,T]$ and hence admits a maximum.
    Inserting the definition for $A_{\mathbb{E}}$, $\Gamma_{\mathbb{E}}$, $R_{\mathbb{E}}$, and $Q_{\mathbb{E}}$ and applying the Cauchy Schwarz inequality yields the assertion.
\end{proof}
Similar estimates can be established for the reconstruction of the state via $\widehat{x}_{\varnothing}$ and $\widehat{x}_{\mathrm{E}}$ as presented in \Cref{subs: expKal} and \Cref{subs: MinEner}, respectively. For this purpose, let $k \in \{ 1,\dots,N \}$. Applying the results of \Cref{subs: TechPrep} with 
\begin{equation*}
    A_* = A_k,
    ~~~
    \Gamma_*= \Gamma_k,
    ~~~
    R_*= R_k,
    ~~~
    Q_* = Q_k,
\end{equation*}
we find that \Cref{lem: est_Kal} yields the existence of $c_{x_k} > 0$ such that for any $t \in [0,T]$ it holds
\begin{equation}\label{eq: est_xk}
    \Vert \widehat{x}_k(t) - \widehat{x}_{\bar{\sigma}} (t) \Vert
    \leq 
    c_{x_k}
    \left(
        \Vert A_k - A_{\bar{\sigma}} \Vert
        +\Vert \Gamma_k - \Gamma_{\bar{\sigma}} \Vert
        +\Vert R_k - R_{\bar{\sigma}} \Vert
        +\Vert Q_k - Q_{\bar{\sigma}} \Vert
    \right).
\end{equation}
The estimate for $\widehat{x}_{\varnothing}$ readily follows.
\begin{proposition}\label{prop: est_varno}
    There exists a constant $c_\varnothing>0$ independent of $N$ such that for all $t \in [0,T]$ we obtain that
    \begin{equation*}
    \begin{aligned}
        \Vert \widehat{x}_{\varnothing}(t) - \widehat{x}_{\bar{\sigma}}(t) \Vert_{P_{\bar{\sigma}(t)}}
        &\leq 
        c_\varnothing~
        \mathbb{E} \left[
        \Vert \mathcal{S}_\sigma - \mathcal{S}_{\bar{\sigma}} \Vert_1
        \right].
    \end{aligned}
    \end{equation*}
\end{proposition}
\begin{proof}
    Looking at the definition of $\widehat{x}_\varnothing$ in \eqref{eq: def_varnothing} and utilizing \eqref{eq: est_xk} we find
    \begin{equation*}
    \begin{aligned}
        &\Vert \widehat{x}_{\varnothing}(t) - \widehat{x}_{\bar{\sigma}}(t) \Vert_{P_{\bar{\sigma}}(t)}
        \leq
        \Vert P_{\bar{\sigma}}^{\tfrac{1}{2}}(t) \Vert
        ~
        \Vert \frac{1}{N} \sum_{k=1}^N \widehat{x}_k(t) - \widehat{x}_{\bar{\sigma}} (t) \Vert\\
        &\leq
        \Vert P_{\bar{\sigma}}^{\tfrac{1}{2}}(t) \Vert
        ~
        \frac{1}{N} \sum_{k=1}^N
        c_{x_k}
        \left(
        \Vert A_\sigma - A_{\bar{\sigma}} \Vert
        +\Vert \Gamma_\sigma - \Gamma_{\bar{\sigma}} \Vert
        +\Vert R_\sigma - R_{\bar{\sigma}} \Vert
        +\Vert Q_\sigma - Q_{\bar{\sigma}} \Vert
        \right).
    \end{aligned}
    \end{equation*}
    The assertion follows with 
    $c_\varnothing 
    = \left( \max\limits_{t \in [0,T]}  \Vert P_{\bar{\sigma}}^{\tfrac{1}{2}}(t) \Vert \right)
    \left( \max\limits_{k \in \{ 1,\dots,N \}} c_{x_k} \right)$. 
\end{proof}
Before we can establish an analogous estimate for the minimizer of the expected energy $\widehat{x}_{\mathrm{E}}$ discussed in \Cref{subs: MinEner} we take a closer look at the precision matrices $P_k(t)$. In the following, let $\mathcal{P}(\cdot) =\sum_{k=1}^N P_k(t)$ be as in \Cref{prop: char_en_min}.
\begin{lemma}\label{lem: estMathcP}
    There exists a constant $c_{\mathcal{P}} > 0$ independent of $N$ such that for all $t \in [0,T]$ it holds
    \begin{equation*}
        \Vert \mathcal{P}^{-1}(t) \Vert 
        \leq 
        \frac{c_{\mathcal{P}}}{N}.
    \end{equation*}
\end{lemma}
\begin{proof}
    For all non-trivial $x \in \mathbb{R}^n$ and $t \in [0,T]$ we find
    \begin{equation*}
        \left( x, \mathcal{P}(t) x \right) 
        =
        \sum_{k=1}^N \left( x,P_k(t)x \right)
        \geq \sum_{k=1}^N \min_{s \in [0,T]} \left( x,P_k(s)x \right),
    \end{equation*}
    where the existence of the minima is ensured by the continuity of the matrices. We denote the (possibly non-unique) time at which the minimum of the $k$-th summand is assumed by $s_k^*$. In addition we denote by $\lambda_k > 0$ the smallest eigenvalue of $P_k(s_k^*)$, where the positivity follows from the positive definiteness of the precision matrices. Denoting $\bar{\lambda} = \max\limits_{k=1,\dots,N} \lambda_k > 0$ it follows
    \begin{equation*}
        \left( x, \mathcal{P}(t) x \right)
        \geq
        \sum_{k=1}^N \lambda_k \Vert x \Vert^2
        \geq 
        N \bar{\lambda} \Vert x \Vert^2.
    \end{equation*}
    The proof is completed with an application of the Lax-Milgram Lemma \cite[Thm.~1.1.3 {\&} Rem.~1.1.3]{Cia78}.
\end{proof}
We conclude this section with the announced estimate for $\widehat{x}_{\mathrm{E}}$.
\begin{proposition}\label{prop: est_en}
    There exists a constant $c_\mathrm{E}>0$ independent of $N$ such that for all $t \in [0,T]$ it holds
    \begin{equation*}
    \begin{aligned}
        \Vert \widehat{x}_{\mathrm{E}}(t) - \widehat{x}_{\bar{\sigma}}(t) \Vert_{P_{\bar{\sigma}(t)}}
        \leq
        c_\mathrm{E}~
        \mathbb{E} \left[
        \Vert \mathcal{S}_\sigma - \mathcal{S}_{\bar{\sigma}} \Vert_1
        \right].
    \end{aligned}
    \end{equation*}
\end{proposition}
\begin{proof}
    With \Cref{prop: char_en_min} for every $t \in [0,T]$ we obtain
    \begin{equation*}
        \Vert \widehat{x}_{\mathrm{E}}(t) - \widehat{x}_{\bar{\sigma}}(t) \Vert_{P_{\bar{\sigma}(t)}}
        \leq
        \Vert P_{\bar{\sigma}}^{\tfrac{1}{2}}(t) \Vert
        \left\Vert \mathcal{P}^{-1}(t) \right\Vert
        \sum_{k=1}^N \Vert P_k(t) \Vert \left\Vert \widehat{x}_k(t) - \widehat{x}_{\bar{\sigma}}(t) \right\Vert.
    \end{equation*}
    The continuity of $P_k(\cdot)$ implies existence of $M>0$ such that for all $t \in [0,T]$ and $k = 1,\dots,N$ it holds
    $\Vert P_k(t) \Vert \leq M$.
    Utilizing \eqref{eq: est_xk} and \Cref{lem: estMathcP} yields
    \begin{equation*}
    \begin{aligned}
        &\Vert \widehat{x}_{\mathrm{E}}(t) - \widehat{x}_{\bar{\sigma}}(t) \Vert_{P_{\bar{\sigma}(t)}}\\
        &\leq
        M \frac{c_{\mathcal{P}}}{N}
        \Vert P_{\bar{\sigma}}^{\tfrac{1}{2}}(t) \Vert
        \sum_{k=1}^N 
        \left( 
        \Vert A_k - A_{\bar{\sigma}} \Vert
        +\Vert \Gamma_k - \Gamma_{\bar{\sigma}} \Vert
        +\Vert R_k - R_{\bar{\sigma}} \Vert
        +\Vert Q_k - Q_{\bar{\sigma}} \Vert
        \right)
    \end{aligned}
    \end{equation*}
    and the claim follows with $c_{\mathrm{E}} = M c_{\mathcal{P}} \max\limits_{s\in[0,T]} \Vert P_{\bar{\sigma}}^{\tfrac{1}{2}}(s) \Vert$.
\end{proof}
The worst case error estimates presented in \Cref{prop: est_exp}, \Cref{prop: est_varno}, and \Cref{prop: est_en} immediately imply bounds for the expected errors. 
\begin{corollary}
    Let $\widehat{x}_{\mathbb{E}}$, $\widehat{x}_\varnothing$, and $\widehat{x}_\mathrm{E}$ be the estimators introduced in \Cref{sec: KF for uncertain systems} and let $c_\mathbb{E}$, $c_\varnothing$, and $c_{\mathrm{E}}$ be the constants from \Cref{prop: est_exp}, \Cref{prop: est_varno}, and \Cref{prop: est_en}, respectively. For all $t \in [0,T]$ it holds
    \begin{equation*}
    \begin{aligned}
        \mathbb{E} \left[
        \Vert \widehat{x}_{\mathbb{E}}(t) - \widehat{x}_{{\sigma}}(t) \Vert_{P_{{\sigma}(t)}}
        \right]
        &\leq
        c_\mathbb{E}~
        \mathbb{E}_{\bar{\sigma}} \left[
        \mathbb{E}_{\sigma} \left[
        \Vert \mathcal{S}_\sigma - \mathcal{S}_{\bar{\sigma}} \Vert_1 \right] \right],\\
        \mathbb{E} \left[
        \Vert \widehat{x}_{\varnothing}(t) - \widehat{x}_{{\sigma}}(t) \Vert_{P_{{\sigma}(t)}}
        \right] 
        &\leq
        c_\varnothing~
        \mathbb{E}_{\bar{\sigma}} \left[
        \mathbb{E}_{\sigma} \left[
        \Vert \mathcal{S}_\sigma - \mathcal{S}_{\bar{\sigma}} \Vert_1 \right] \right],\\
        \mathbb{E} \left[
        \Vert \widehat{x}_{\mathrm{E}}(t) - \widehat{x}_{{\sigma}}(t) \Vert_{P_{{\sigma}(t)}}
        \right]
        &\leq
        c_\mathrm{E}~
        \mathbb{E}_{\bar{\sigma}} \left[
        \mathbb{E}_{\sigma} \left[
        \Vert \mathcal{S}_\sigma - \mathcal{S}_{\bar{\sigma}} \Vert_1 \right] \right].
    \end{aligned}
    \end{equation*}
\end{corollary}
\begin{proof}
    The proof follows by forming the expectation on both sides of the corresponding worst case error estimates. 
\end{proof}
We found that all three state estimators discussed in this work satisfy the same qualitative worst case error bound. They express the fact that these errors can be quantified by the level of variance of the parameter-dependent matrices $A_{\bar \sigma},\Gamma_{\bar \sigma},R_{\bar \sigma},Q_{\bar \sigma}$ .

%
\section{Numerical experiments}\label{sec: numExp}
%
In this section we numerically realize and compare the three estimators for two examples of uncertain linear dynamical systems.
%
\subsection{General setup}\label{subs: genSetup}
%
We consider disturbed uncertain linear systems of the form 
\begin{equation}\label{eq: numSys}
\begin{alignedat}{2}
    \dot{x}(t) &= A_\sigma x(t) + B v(t) 
    ~~~~~&& t \in (0,T),\\
    x(0) &= x_0 + \eta, &&\\
    y(t) &= C x(t) + \mu(t)
    ~~~~~&& t \in (0,T),
\end{alignedat}
\end{equation}
where $\sigma \in \Sigma_A \subset \mathbb{R}^{s_A}$ with $ \vert \Sigma_A \vert = N_A$, as in \Cref{subsec: ModUnc}. We associate positive definite weighting/covariance matrices $\Gamma$, $R$, and $Q$ with the three unknown errors $\eta$, $v$, and $\mu$. Here we do not investigate uncertainties in the noise covariance matrices 
implying $N_\Gamma = N_R = N_Q = 1$ leading to $N = N_A$.

All ODEs are solved using an equidistant grid $0 = t_0 < t_1 < \dots < t_{1000} = T$ and the \matlab\, solver \texttt{ode15s} with a relative tolerance of $10^{-8}$. In particular, Kalman filter equations of the form \eqref{eq: KF_fixed}-\eqref{eq: Ricc_fixed} are treated as $n+n^2$ dimensional ODEs. For systems of higher dimension $n$ this can certainly be improved upon by taking advantage of the decoupling of the two equations and utilizing solvers tailored to Riccati equations. 

The measured output $y \in \mathcal{L}_0^T$ is artificially generated as follows. Using the \matlab\, function \texttt{normrnd} we construct realizations of $\eta \sim \mathcal{N}(0,\Gamma)$, $v(t_k) \sim \mathcal{N}(0,R)$, and $\mu(t_k) \sim \mathcal{N}(0,Q)$ for $k = 0,\dots,1000$ from which $v \in \mathcal{L}_0^T$, $\mu \in \mathcal{L}_0^T$ are obtained via linear interpolation. Finally, we designate some $\bar{\sigma} \in \Sigma_A$ to be the true, unknown parameter and obtain $y$ by solving \eqref{eq: numSys} based on the constructed errors and the parameter $\bar{\sigma}$. 
\begin{remark}
    We point out that this realization of the unknown errors based on normal distributions can only be done rigorously in discrete time. The extension to the continuous time setting via interpolation is of heuristic nature as a construction in the spirit of $ v(t) \sim \mathcal{N}(0,R) $, $t \in (0,T)$ leads to the stochastic formulation and the regularity gap discussed in \Cref{rem: L2gap}.
\end{remark}

The construction of the estimators requires the realization of the $N_A$ individual Kalman filters $ \widehat{x}_k $. These are independent of each other, and therefore they are computed in parallel using \texttt{parfor} from the \matlab\, Parallel Computing Toolbox.

The results of this section were realized in \matlab\, R2024b and the computations were executed on a Lenovo ThinkPad E14 Intel(R) Core(TM) Ultra 7 155H 1.40 GHz with 32 GB memory. The code is available in \cite{Sch25}.
%
\subsection{Harmonic oscillator}
%
As a first example we consider a disturbed harmonic oscillator modeled as
\begin{equation*}
\begin{aligned}
    m \ddot{z}(t) + \sigma \dot{z}(t) + k z(t) &= v(t),
    ~~~~~~~~~
    t \in (0,T)\\
    z(0) = x_{0,1} + \eta_1,~ \dot{z}(0) &= x_{0,2} + \eta_2, 
\end{aligned}
\end{equation*}
where $z(t)$ represents the position at time $t$ and $m>0$, $\sigma \geq 0$, and $k \geq 0$ represent the mass, damping coefficient, and spring constant, respectively. The function $v \in \mathcal{L}_0^T$ represents the unknown disturbance in the dynamics. The modeled initial position and velocity are given by $x_{0,1}$ and $x_{0,2}$ and are subject to the initial errors $\eta_1$ and $\eta_2$. Further, we assume access to a disturbed measurement of the position modeled as 
\begin{equation*}
    y(t) = z(t) + \mu(t) 
    ~~~ t \in (0,T),
\end{equation*}
with a deterministic but unknown output error $\mu \in \mathcal{L}_0^T$. In first order form the system reads
\begin{equation*}
\begin{alignedat}{2}
    \dot{x}(t) &= 
    \begin{bmatrix}
        0 & 1 \\
        - \frac{k}{m} & - \frac{\sigma}{m}
    \end{bmatrix}
    x(t)
    +
    \begin{bmatrix} 0 \\ 1 \end{bmatrix} v(t)
    ~~~ &&t \in (0,T),\\
    x(0) &= x_0 + \eta,&&\\ 
    y(t) & = \begin{bmatrix} 1 & 0 \end{bmatrix} x(t)
    + \mu(t)
    ~~~ &&t \in (0,T),
\end{alignedat}
\end{equation*}
with modeled initial state $x_0 = \begin{bmatrix} x_{0,1} & x_{0,2} \end{bmatrix}^\top$.

\begin{figure}
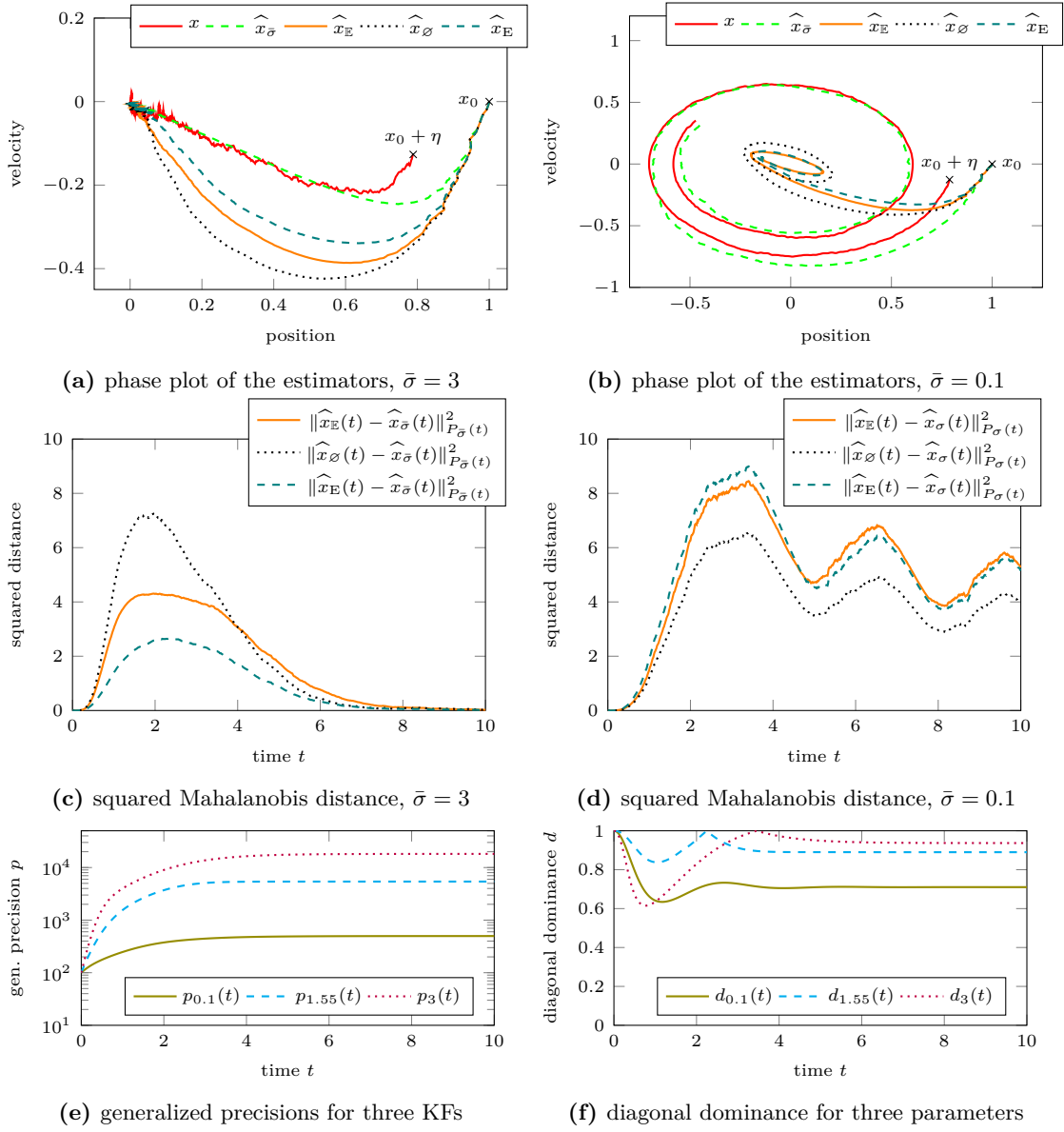

    \begin{subfigure}{0.49\textwidth}
		\begin{tikzpicture}
            \begin{axis}[
            width=\textwidth,
            height=0.25 \textheight,
            xlabel={position},
            ylabel={velocity},
            legend style={  at={(1.05,1.05)},
                            anchor=north east,
                            legend columns=5,
                            font=\tiny
            },
            grid=minor,
            xmin = -0.1,
            xmax = 1.05,
            ymin = -0.45,
            ymax = 0.2,
            label style={font=\tiny},
            tick label style={font=\tiny},
            ]
            \input{Data/Oscillator/sigma3/x_tr_pha} 
            \input{Data/Oscillator/sigma3/kf_tr_pha} 
            \input{Data/Oscillator/sigma3/kf_exp_pha} 
            \input{Data/Oscillator/sigma3/kf_varn_pha} 
            \input{Data/Oscillator/sigma3/kf_en_pha} 
            \node[black,left] at (axis cs:1,0){\tiny{$x_0$}};
            \draw plot[mark=x] coordinates{(1,0)};
            \node[black,above] at (axis cs:0.7884,-0.1266){\tiny{$x_0+\eta$}};
            \draw plot[mark=x] coordinates{(0.7884,-0.1266)};
            \addlegendentry{$x$}
            \addlegendentry{$\widehat{x}_{\bar{\sigma}}$}
            \addlegendentry{$\widehat{x}_{\mathbb{E}}$}
            \addlegendentry{$\widehat{x}_{\varnothing}$}
            \addlegendentry{$\widehat{x}_{\mathrm{E}}$}
            \end{axis}
        \end{tikzpicture}
        \caption{phase plot of the estimators, $\bar{\sigma} = 3$}
        \label{subf: pos3}
	\end{subfigure}
    \begin{subfigure}{0.49\textwidth}
		\begin{tikzpicture}
            \begin{axis}[
            width=\textwidth,
            height=0.25 \textheight,
            xlabel={position},
            ylabel={velocity},
            legend style={  at={(1.05,1.05)},
                            anchor=north east,
                            legend columns=5,
                            font=\tiny
            },
            grid=minor,
            xmin = -0.8,
            xmax = 1.25,
            ymin = -1,
            ymax = 1.2,
            label style={font=\tiny},
            tick label style={font=\tiny}
            ]
            \input{Data/Oscillator/sigma01/x_tr_pha} 
            \input{Data/Oscillator/sigma01/kf_tr_pha} 
            \input{Data/Oscillator/sigma01/kf_exp_pha} 
            \input{Data/Oscillator/sigma01/kf_varn_pha} 
            \input{Data/Oscillator/sigma01/kf_en_pha} 
            \node[black,right] at (axis cs:1,0){\tiny{$x_0$}};
            \draw plot[mark=x] coordinates{(1,0)};
            \node[black,above] at (axis cs:0.7884,-0.1266){\tiny{$x_0+\eta$}};
            \draw plot[mark=x] coordinates{(0.7884,-0.1266)};
            \addlegendentry{$x$}
            \addlegendentry{$\widehat{x}_{\bar{\sigma}}$}
            \addlegendentry{$\widehat{x}_{\mathbb{E}}$}
            \addlegendentry{$\widehat{x}_{\varnothing}$}
            \addlegendentry{$\widehat{x}_{\mathrm{E}}$}
            \end{axis}
        \end{tikzpicture}
        \caption{phase plot of the estimators, $\bar{\sigma} = 0.1$}
        \label{subf: pos01}
	\end{subfigure}
    \begin{subfigure}{0.49\textwidth}
		\begin{tikzpicture}
            \begin{axis}[
            width=\textwidth,
            height=0.25 \textheight,
            xlabel={time $t$},
            ylabel={squared distance},
            legend style={at={(1.05,1.15)}, anchor=north east},
            grid=minor,
            xmin = 0,
            xmax = 10,
            ymin = 0,
            ymax = 10,
            label style={font=\tiny},
            tick label style={font=\tiny},
            legend style={font=\tiny}
            ]
            \input{Data/Oscillator/sigma3/mah_dist_exp}
            \input{Data/Oscillator/sigma3/mah_dist_varn}
            \input{Data/Oscillator/sigma3/mah_dist_en} 
            \addlegendentry{$ \Vert \widehat{x}_{\mathbb{E}}(t) -\widehat{x}_{\bar{\sigma}}(t) \Vert_{P_{\bar{\sigma}}(t)}^2$}
            \addlegendentry{$ \Vert \widehat{x}_{\varnothing}(t) -\widehat{x}_{\bar{\sigma}}(t) \Vert_{P_{\bar{\sigma}}(t)}^2 $}
            \addlegendentry{$ \Vert \widehat{x}_{\mathrm{E}}(t) -\widehat{x}_{\bar{\sigma}}(t) \Vert_{P_{\bar{\sigma}}(t)}^2 $}
            \end{axis}
        \end{tikzpicture}
        \caption{squared Mahalanobis distance, $\bar{\sigma} = 3$}
        \label{subf: MahD3}
	\end{subfigure}
    \begin{subfigure}{0.49\textwidth}
		\begin{tikzpicture}
            \begin{axis}[
            width=\textwidth,
            height=0.25 \textheight,
            xlabel={time $t$},
            ylabel={squared distance},
            legend style={at={(1.05,1.15)}, anchor=north east},
            grid=minor,
            xmin = 0,
            xmax = 10,
            ymin = 0,
            ymax = 10,
            label style={font=\tiny},
            tick label style={font=\tiny},
            legend style={font=\tiny}
            ]
            \input{Data/Oscillator/sigma01/mah_dist_exp}
            \input{Data/Oscillator/sigma01/mah_dist_varn}
            \input{Data/Oscillator/sigma01/mah_dist_en} 
            \addlegendentry{$ \Vert \widehat{x}_{\mathbb{E}}(t) -\widehat{x}_{\sigma}(t) \Vert_{P_{\sigma}(t)}^2$}
            \addlegendentry{$ \Vert \widehat{x}_{\varnothing}(t) -\widehat{x}_{\sigma}(t) \Vert_{P_{\sigma}(t)}^2 $}
            \addlegendentry{$ \Vert \widehat{x}_{\mathrm{E}}(t) -\widehat{x}_{\sigma}(t) \Vert_{P_{\sigma}(t)}^2 $}
            \end{axis}
        \end{tikzpicture}
        \caption{squared Mahalanobis distance, $\bar{\sigma} = 0.1$}
        \label{subf: Mah01}
	\end{subfigure}
    \begin{subfigure}{0.49\textwidth}
		\begin{tikzpicture}
            \begin{semilogyaxis}[
            width=\textwidth,
            height=0.2 \textheight,
            xlabel={time $t$},
            ylabel={gen.~precision $p$},
            legend style={  at={(0.95,0.05)},
                            anchor=south east,
                            legend columns=3,
                            font=\tiny
            },
            grid=minor,
            grid style={white},
            xmin = 0,
            xmax = 10,
            ymin = 10,
            ymax = 5*10^4,
            label style={font=\tiny},
            tick label style={font=\tiny}
            ]
            \input{Data/Oscillator/gen_prec_01}
            \input{Data/Oscillator/gen_prec_155}
            \input{Data/Oscillator/gen_prec_3} 
            \addlegendentry{$p_{0.1}(t)$}
            \addlegendentry{$p_{1.55}(t)$}
            \addlegendentry{$p_{3}(t)$}
            \end{semilogyaxis}
        \end{tikzpicture}
        \caption{generalized precisions for three KFs}
        \label{subf: genPrec}
	\end{subfigure}
    \begin{subfigure}{0.49\textwidth}
		\begin{tikzpicture}
            \begin{axis}[
            width=\textwidth,
            height=0.2 \textheight,
            xlabel={time $t$},
            ylabel={diagonal dominance $d$},
            legend style={  at={(0.95,0.05)},
                            anchor=south east,
                            legend columns=3,
                            font=\tiny
            },
            grid=minor,
            grid style={white},
            xmin = 0,
            xmax = 10,
            ymin = 0,
            ymax = 1,
            label style={font=\tiny},
            tick label style={font=\tiny}
            ]
            \input{Data/Oscillator/diagDom_01}
            \input{Data/Oscillator/diagDom_155}
            \input{Data/Oscillator/diagDom_3} 
            \addlegendentry{$d_{0.1}(t)$}
            \addlegendentry{$d_{1.55}(t)$}
            \addlegendentry{$d_{3}(t)$}
            \end{axis}
        \end{tikzpicture}
        \caption{diagonal dominance for three parameters}
        \label{subf: diagDom}
	\end{subfigure}
    \caption{Harmonic oscillator with uncertain damping}
    \label{fig: oscillator}
\end{figure}

For our experiments we fix the mass and spring constant to $m = k = 1$ and set the time horizon to $T = 10$. As an undisturbed initial state we set $x_0 = \begin{bmatrix} 1 & 0 \end{bmatrix}^\top$ and the covariances of the three errors are given as $\Gamma = 0.1 \, \mathrm{Id} $, $R = 0.05$, and $Q = 0.05$. The uncertainty of the dynamics lies in the damping parameter $\sigma$ and we set $N = N_A = 101$ and $\Sigma_A = \{ 0.1 + \tfrac{k}{N_A-1} 2.9 \colon k = 0,\dots,N_A-1 \} $. We generate two outputs based on the parameters $\bar{\sigma} = 3$ and $\bar{\sigma} = 0.1$ as described in the previous subsection. The results are presented in \Cref{fig: oscillator} and we observe the following.

We first consider the setting of $\bar{\sigma} = 3$, i.e., the underlying true damping parameter is the largest appearing in the family $\Sigma_A$. In \Cref{subf: pos3} we present the phase plots of the three estimators $\widehat{x}_{\mathbb{E}}$, $\widehat{x}_{\varnothing}$, and $\widehat{x}_\mathrm{E}$ resulting from the associated output $y$ and compare them with the Kalman filter $\widehat{x}_{\bar{\sigma}}$ constructed based on the hidden parameter $\bar{\sigma} = 3$ and the associated state trajectory $x$.
The trajectories starting in their respective initial states marked as crosses stabilize in a region around the origin.
As expected, the Kalman filter $\widehat{x}_{\bar{\sigma}}$ successfully approximates the state trajectory $x$. The plot further illustrates that among the estimators constructed without knowledge of the hidden parameter the energy minimizer $\widehat{x}_\mathrm{E}$ yields the best approximation of $\widehat{x}_{\bar{\sigma}}$ with respect to the Euclidean distance. In addition \Cref{subf: MahD3} shows that the energy minimizer outperforms the two remaining estimators with respect to the Mahalanobis distance. 
We proceed with the results associated with $\bar{\sigma} = 0.1$. Here the true damping parameter is the smallest appearing in $\Sigma_A$. The associated phase plot presented in \Cref{subf: pos01} shows that all three estimators $\widehat{x}_{\mathbb{E}}$, $\widehat{x}_{\varnothing}$, and $\widehat{x}_\mathrm{E}$ noticeably deviate from the Kalman filter $\widehat{x}_{\bar{\sigma}}$. The associated Mahalanobis differences are displayed in \Cref{subf: Mah01} and we find that, in contrast to the case of $\bar{\sigma} = 3$, 
here the energy minimizer $\widehat{x}_\mathrm{E}$ does not outperform $\widehat{x}_\mathbb{E}$ and $\widehat{x}_\varnothing$. Note, however, that the accuracies of the three are noticeably more similar than in the previous setting. 
In an attempt to interpret these results the reader is reminded of the formulas for $\widehat{x}_\varnothing$ and $\widehat{x}_\mathrm{E}$. According to \eqref{eq: def_varnothing} the estimator $\widehat{x}_\varnothing$ is given as an (unweighted) mean of the family of Kalman filters $x_k$. On the other hand the energy minimizer $\widehat{x}_\mathrm{E}$ is characterized as a weighted mean where the precision matrices $P_k$ act as weights. In fact, if all $P_k$ were diagonal, each component of $\widehat{x}_\mathrm{E}$ would be a mean of the associated components of $\widehat{x}_k$ weighted by their respective precisions. In such a setting one expects $\widehat{x}_\mathrm{E}$ to strongly resemble the family members $\widehat{x}_k$ associated with high precisions, while family members with lower precisions are neglected. 
While the precision matrices can not be expected to be diagonal, we do find that for this example they are diagonally dominant. To illustrate this, for a given parameter $\sigma \in \Sigma_A$ we consider the quantity
\begin{equation*}
    0 \leq 
    d_\sigma(t) = \min_{j = 1,\dots,n} \frac{\vert p_{jj}(t) \vert}{\sum_{i=1}^n \vert p_{ji}(t) \vert}
    \leq 1,
\end{equation*}
where $p{ij}(t)$ are the entries of $P_\sigma(t)$. Note that since the precision matrices are positive definite, the denominator is non zero. Now $P_\sigma(t)$ is diagonally dominant if $d_\sigma(t) \geq 0.5$ and diagonal if $d_\sigma(t) = 1$.
The announced diagonal dominance of the precision matrices is illustrated in \Cref{subf: diagDom}. 
In \Cref{subf: genPrec} we display the generalized precisions (cf., \Cref{def: concStoch})
\begin{equation*}
    p_\sigma(t) = \det(P_\sigma(t))
\end{equation*}
to quantify the precisions associated with the family members and find that it depends monotonously increasing on the damping parameter. This fits the intuition that larger damping parameters lead to stronger dissipation of errors. Note also that the precisions increase over time reflecting the fact that the effect of the initial error $\eta $ decreases over time. 

These considerations offer an explanation for the observations made while comparing the effects of the hidden parameters $\bar{\sigma} = 3$ and $\bar{\sigma} = 0.1$. Since the precision matrices are diagonally dominant, we expect the estimators to behave similarly to the setting of diagonal precision matrices. Hence the weighted mean $\widehat{x}_\mathrm{E}$ favors family members with higher precisions, i.e., trajectories associated with larger damping parameters. Finally, we note that all three designed estimators display noticeable differences from the Kalman filter based on the hidden parameter. This, however, is not surprising as the set of damping parameters leads to a family of system matrices $A_\sigma$ representing a wide range of dynamics, making the state reconstruction a rather challenging task. 

%
\subsection{Connected amplidynes}
%

\begin{figure}
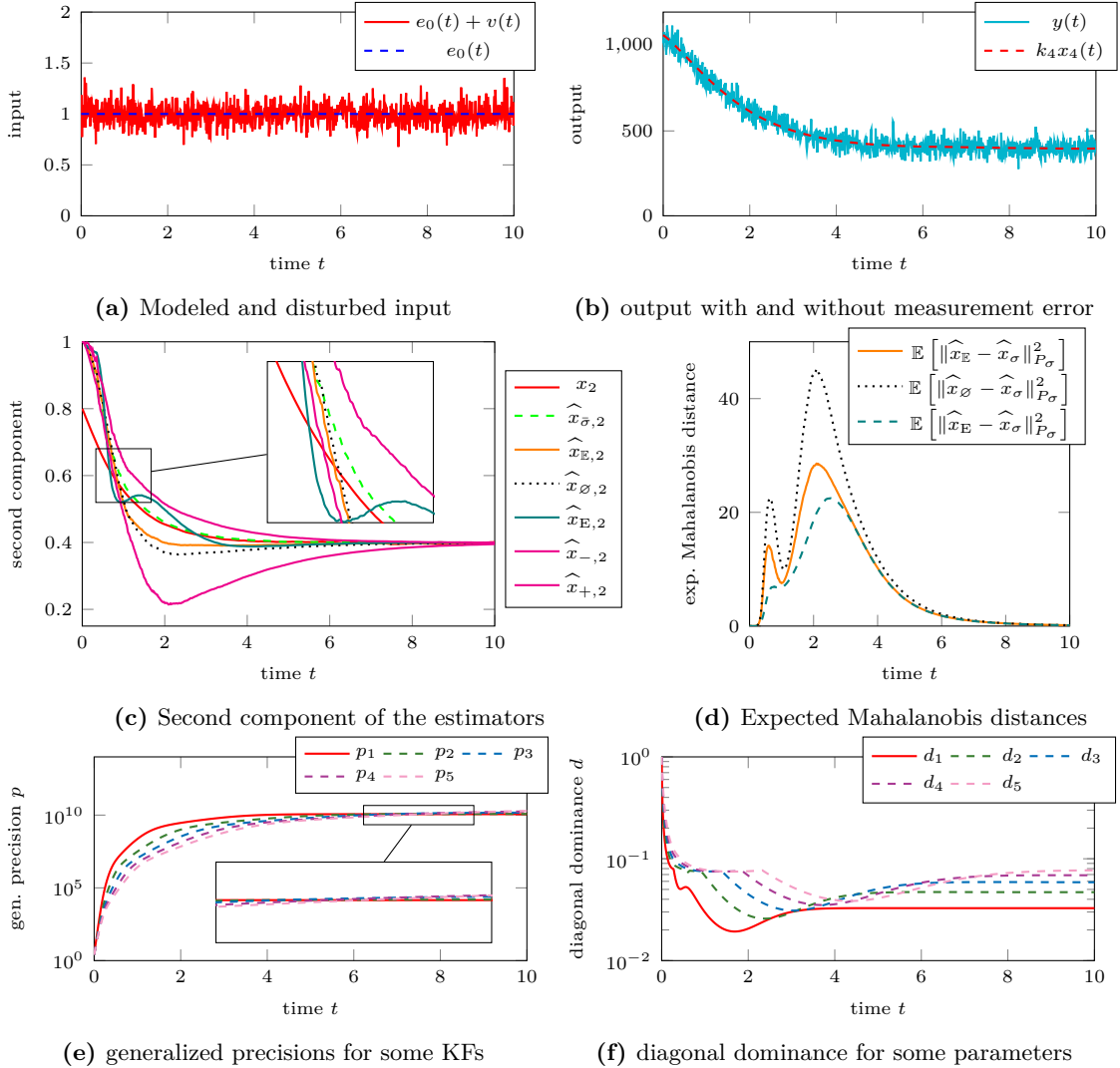

    \begin{subfigure}{0.49\textwidth}
		\begin{tikzpicture}
            \begin{axis}[
            width=\textwidth,
            height=0.2 \textheight,
            xlabel={time $t$},
            ylabel={input},
            legend style={  at={(1.05,1.05)},
                anchor=north east,
                legend columns=1,
                font=\tiny
            },
            grid=minor,
            grid style={white},
            xmin = 0,
            xmax = 10,
            ymin = 0,
            ymax = 2,
            label style={font=\tiny},
            tick label style={font=\tiny}
            ]
            \input{Data/Amplidyne/input_dist}
            \input{Data/Amplidyne/intput_mod}
            \addlegendentry{$e_0(t)+v(t)$}
            \addlegendentry{$e_0(t)$}
            \end{axis}
        \end{tikzpicture}
        \caption{Modeled and disturbed input}
        \label{subf: Input}
	\end{subfigure}
    \begin{subfigure}{0.49\textwidth}
		\begin{tikzpicture}
            \begin{axis}[
            width=\textwidth,
            height=0.2 \textheight,
            xlabel={time $t$},
            ylabel={output},
            legend style={  at={(1.05,1.05)},
                anchor=north east,
                legend columns=1,
                font=\tiny
            },
            grid=minor,
            grid style={white},
            xmin = 0,
            xmax = 10,
            ymin = 0,
            ymax = 1200,
            label style={font=\tiny},
            tick label style={font=\tiny}
            ]
            \input{Data/Amplidyne/output_dist}
            \input{Data/Amplidyne/output_tr}
            \addlegendentry{$y(t)$}
            \addlegendentry{$k_4 x_4(t)$}
            \end{axis}
        \end{tikzpicture}
        \caption{output with and without measurement error}
        \label{subf: Output}
	\end{subfigure}
    \begin{subfigure}{0.59\textwidth}
		\begin{tikzpicture}[spy using outlines={magnification=2.5, connect spies}]
            \begin{axis}[
            width=0.8\textwidth,
            height=0.25 \textheight,
            xlabel={time $t$},
            ylabel={second component},
            legend style={  at={(1.025,0.9)},
                anchor=north west,
                legend columns=1,
                font=\tiny
            },
            grid=minor,
            grid style={white},
            xmin = 0,
            xmax = 10,
            ymin = 0.15,
            ymax = 1,
            label style={font=\tiny},
            tick label style={font=\tiny}
            ]
            \input{Data/Amplidyne/x_tr}
            \input{Data/Amplidyne/kf_tr}
            \input{Data/Amplidyne/kf_exp}
            \input{Data/Amplidyne/kf_varn}
            \input{Data/Amplidyne/kf_en}
            \input{Data/Amplidyne/kf_min}
            \input{Data/Amplidyne/kf_max}
            \addlegendentry{$x_2$}
            \addlegendentry{$\widehat{x}_{\bar{\sigma},2}$}
            \addlegendentry{$\widehat{x}_{\mathbb{E},2}$}
            \addlegendentry{$\widehat{x}_{\varnothing,2}$}
            \addlegendentry{$\widehat{x}_{\mathrm{E},2}$}
            \addlegendentry{$\widehat{x}_{-,2}$}
            \addlegendentry{$\widehat{x}_{+,2}$}
            %
            \newcommand*\zoomFactor{1.732050807568877} 
            \newcommand*\spyPoint{axis cs:1,0.6}
            \newcommand*\zoomPoint{axis cs:6.5,0.7}
            %
            \node[very thin, rectangle, draw, height=0.0333333333333333\textheight,width=0.0833333333333333\textwidth] (spyrectangle_a) at (\spyPoint) {};
            \node[rectangle, draw, height=0.1\textheight,width=0.25\textwidth] (zoomrectangle_a) at (\zoomPoint) {};
            \draw (spyrectangle_a) edge (zoomrectangle_a);

            \begin{scope}
                \clip (zoomrectangle_a.south west) rectangle (zoomrectangle_a.north east);
                \pgfmathparse{\zoomFactor^2/(\zoomFactor-1)}
                \begin{scope}[scale around={\zoomFactor:($(\zoomPoint)!\zoomFactor^2/(\zoomFactor^2-1)!(\spyPoint)$)}]
                    \input{Data/Amplidyne/x_tr}
                    \input{Data/Amplidyne/kf_tr}
                    \input{Data/Amplidyne/kf_exp}
                    \input{Data/Amplidyne/kf_varn}
                    \input{Data/Amplidyne/kf_en}
                    \input{Data/Amplidyne/kf_min}
                    \input{Data/Amplidyne/kf_max}
                \end{scope}
            \end{scope}
            \end{axis}
        \end{tikzpicture}
        \caption{Second component of the estimators}
        \label{subf: secComp}
	\end{subfigure}
    \begin{subfigure}{0.39\textwidth}
		\begin{tikzpicture}
            \begin{axis}[
            width=\textwidth,
            height=0.25 \textheight,
            xlabel={time $t$},
            ylabel={exp.~Mahalanobis distance},
            legend style={  at={(1.05,1.05)},
                            anchor=north east,
                            legend columns=1,
                            font=\tiny
            },
            grid=minor,
            grid style={white},
            xmin = 0,
            xmax = 10,
            ymin = 0,
            ymax = 50,
            label style={font=\tiny},
            tick label style={font=\tiny}
            ]
            \input{Data/Amplidyne/exp_Mah_exp}
            \input{Data/Amplidyne/exp_Mah_varn}
            \input{Data/Amplidyne/exp_Mah_en}
            \addlegendentry{$\mathbb{E} \left[ \Vert \widehat{x}_{\mathbb{E}} -\widehat{x}_{\sigma} \Vert_{P_{\sigma}}^2 \right]$}
            \addlegendentry{$\mathbb{E} \left[ \Vert \widehat{x}_{\varnothing} -\widehat{x}_{\sigma} \Vert_{P_{\sigma}}^2 \right]$}
            \addlegendentry{$\mathbb{E} \left[ \Vert \widehat{x}_{\mathrm{E}} -\widehat{x}_{\sigma} \Vert_{P_{\sigma}}^2 \right]$}
            \end{axis}
        \end{tikzpicture}
        \caption{Expected Mahalanobis distances}
        \label{subf: expMah}
	\end{subfigure}
    \begin{subfigure}{0.49\textwidth}
		\begin{tikzpicture}[spy using outlines={magnification=1, connect spies}]
            \begin{semilogyaxis}[
            width=\textwidth,
            height=0.2 \textheight,
            xlabel={time $t$},
            ylabel={gen.~precision $p$},
            legend style={  at={(1.05,1.1)},
                            anchor=north east,
                            legend columns=3,
                            font=\tiny
            },
            grid=minor,
            grid style={white},
            xmin = 0,
            xmax = 10,
            ymin = 1,
            ymax = 10^14,
            label style={font=\tiny},
            tick label style={font=\tiny}
            ]
            \input{Data/Amplidyne/gen_prec_1}
            \input{Data/Amplidyne/gen_prec_2}
            \input{Data/Amplidyne/gen_prec_3}
            \input{Data/Amplidyne/gen_prec_4}
            \input{Data/Amplidyne/gen_prec_5}
            \addlegendentry{$p_1$}
            \addlegendentry{$p_2$}
            \addlegendentry{$p_3$}
            \addlegendentry{$p_4$}
            \addlegendentry{$p_5$}
            %
            \newcommand*\zoomFactor{1.58113883008419} 
            \newcommand*\spyPoint{axis cs:7.5,10^10}
            \newcommand*\zoomPoint{axis cs:6,10^4}
            %
            \node[very thin, rectangle, draw, height=height=0.02\textheight,width=0.2\textwidth] (spyrectangle_b) at (\spyPoint) {};
            \node[rectangle, draw, height=0.05\textheight,width=0.5\textwidth] (zoomrectangle_b) at (\zoomPoint) {};
            \draw (spyrectangle_b) edge (zoomrectangle_b);
            
            \begin{scope}
                \clip (zoomrectangle_b.south west) rectangle (zoomrectangle_b.north east);
                \begin{scope}[scale around={1.58113883008419:(axis cs:8.5,10^14)}]
                    \input{Data/Amplidyne/gen_prec_1}
                    \input{Data/Amplidyne/gen_prec_2}
                    \input{Data/Amplidyne/gen_prec_3}
                    \input{Data/Amplidyne/gen_prec_4}
                    \input{Data/Amplidyne/gen_prec_5}
                \end{scope}
            \end{scope}
            \end{semilogyaxis}
        \end{tikzpicture}
        \caption{generalized precisions for some KFs}
        \label{subf: genPrecA}
	\end{subfigure}
    \begin{subfigure}{0.49\textwidth}
		\begin{tikzpicture}
            \begin{semilogyaxis}[
            width=\textwidth,
            height=0.2 \textheight,
            xlabel={time $t$},
            ylabel={diagonal dominance $d$},
            legend style={  at={(1.05,1.1)},
                            anchor= north east,
                            legend columns=3,
                            font=\tiny
            },
            grid=minor,
            grid style={white},
            xmin = 0,
            xmax = 10,
            ymin = 0.01,
            ymax = 1,
            label style={font=\tiny},
            tick label style={font=\tiny}
            ]
            \input{Data/Amplidyne/diagDom_1}
            \input{Data/Amplidyne/diagDom_2}
            \input{Data/Amplidyne/diagDom_3}
            \input{Data/Amplidyne/diagDom_4}
            \input{Data/Amplidyne/diagDom_5}
            \addlegendentry{$d_1$}
            \addlegendentry{$d_2$}
            \addlegendentry{$d_3$}
            \addlegendentry{$d_4$}
            \addlegendentry{$d_5$}
            \end{semilogyaxis}
        \end{tikzpicture}
        \caption{diagonal dominance for some parameters}
        \label{subf: diagDomA}
	\end{subfigure}
    \caption{Connected amplidynes with uncertain inductances}
    \label{fig: amplidyne}
\end{figure}

Our second example is given by an electrical circuit that amplifies a given input. An amplidyne returns an amplified copy of the input signal and can be modeled as a linear dynamical system, see \cite[Ch.~1.12]{KwaSiv72}. We consider a system given by two connected amplidynes, i.e., the output of the first one serves as an input of the second one. Hence the resulting circuit consists of four components and the associated dynamics are modeled as 
\begin{equation*}
\begin{alignedat}{2}
    \dot{x}(t) &= 
    \begin{bmatrix}
        - \frac{\rho_1}{L_1} & 0 & 0 & 0\\
        \frac{k_1}{L_2} & - \frac{\rho_2}{L_2} & 0 & 0 \\
        0 & \frac{k_2}{L_3} & - \frac{\rho_3}{L_3} & 0 \\ 
        0 & 0 & \frac{k_3}{L_4} & - \frac{\rho_4}{L_4}
    \end{bmatrix}
    x(t)
    +
    \begin{bmatrix}
        \frac{e_0 (t)}{L_1} \\ 0 \\ 0 \\ 0
    \end{bmatrix}
    +
    \begin{bmatrix}
        \frac{1}{L_1} \\ 0 \\ 0 \\ 0
    \end{bmatrix}
    v(t)
    ~~~~~&& t \in (0,T),\\
    x(0) &=
    x_0 + \eta, && \\
    y(t) 
    &= 
    \begin{bmatrix}
        0 & 0 & 0 & k_4
    \end{bmatrix}
    x(t)
    + \mu(t)
    ~~~~~&& t \in (0,T),
\end{alignedat}
\end{equation*}
where $x(t) \in \mathbb{R}^4$ describes the current of the four components at time $t$. The positive parameters $\rho_i$ and $L_i$ represent the resistance and inductance of the $i$-th component, respectively. Further the constants $k_i > 0$ describe the proportion of the currents $x_i$ and the voltages $e_i$ in the components, i.e., for $i=1,2,3,4$ it holds $e_i = k_i x_i$. The known, time-dependent input entering the first component is denoted by $e_0$ and is subject to the disturbance $v$. Finally, as an output of the system we measure the output of the second amplidyne $e_4 = k_4 x_4$.
We note that strictly speaking our formulation of the Kalman filter does not allow for a known forcing term like the input $e_0$ of the first amplidyne. This, however, can be incorporated without any major challenges.  
For our implementation we fix $\rho_1 = \rho_3 = 5$, $\rho_2 = \rho_4 = 10$, $k_1 = k_3 = 20$, $k_2 = k_4 = 50$, $L_1 = 0.5$, $T=10$, and $e_0(t) \equiv 1$. The undisturbed initial state is set to $x_0 = \begin{bmatrix} 0.5 & 1 & 10 & 20 \end{bmatrix}^\top$ and the uncertainty of the system lies in the inductances $\sigma = \begin{bmatrix} L_2 & L_3 & L_4 \end{bmatrix}^\top \in \mathbb{R}^3$. We define 
$\Sigma_1 = \{ 10, 12.5, 15, 17.5, 20 \}$,
$\Sigma_2 = \{ 0.5, 0.75, 1, 1.25, 1.5 \}$,
and 
$\Sigma_3 = \{ 10, 17.5, 25, 32.5, 40 \}$
and construct $\Sigma_A = \{ \sigma \in \mathbb{R}^3 \colon \sigma_k \in \Sigma_k, k = 1,\dots,3 \} $. We obtain $N = N_A = 125$ and $s_A = 3$. 
The output $y$ is generated as described in \Cref{subs: genSetup} according to the covariance matrices 
\begin{equation*}
\begin{aligned}
    \Gamma &= 0.25~ \mathrm{diag} (\vert x_{0,1} \vert,\vert x_{0,2} \vert,\vert x_{0,3} \vert,\vert x_{0,4} \vert) = \mathrm{diag} (0.125,0.25,2.5,5),\\
    R &= (0.1 \vert e_0(0) \vert)^2 = 0.01,\\
    Q &= (0.1 * 400 \vert e_0(0) \vert)^2 = 1600,
\end{aligned}
\end{equation*}
and the true parameter $\bar{\sigma} =  \begin{bmatrix} 10 & 0.5 & 10 \end{bmatrix}^\top$. We note that the covariance matrices are of different magnitudes to fit the scales of the different components. We illustrate the relationship between the input $e_0$ and the output $k_4 x_4 = e_4$ and the impact of the errors in \Cref{subf: Input} and \Cref{subf: Output}. We observe that the output of the two connected amplidynes associated with $\bar{\sigma}$ stabilizes at a copy of $e_0$ increased by a factor of 400.  

Our focus lies on the investigation of the structural differences to the first example. In \Cref{subf: secComp} we present the second component of the disturbed trajectory and the various estimators. Additionally we plot in magenta color the minima and maxima of the family of Kalman filters denoted by 
\begin{equation*}
    \widehat{x}_{-,2}(t) = \min_{k=1,\dots,N} \widehat{x}_{k,2}(t)
    ~~~~~\text{and}~~~~~
    \widehat{x}_{+,2}(t) = \max_{k=1,\dots,N} \widehat{x}_{k,2}(t),
\end{equation*}
illustrating their convex hull.
As expected the unweighted average $\widehat{x}_\varnothing$ lies in this convex hull. This is different from the energy minimizer $\widehat{x}_\mathrm{en}$ which 
for some times $t$ has values outside this hull. 
This indicates that the off diagonal entries of the precision matrices $P_k$ play an important role. 
This corresponds to the fact that
the precision matrices are by no means diagonal dominant. In \Cref{subf: diagDomA} we present the quantities 
\begin{equation*} 
    d_k = d_{\sigma_k}(t) = \min_{j = 1,\dots,n} \frac{\vert p_{jj}^{\sigma_k}(t) \vert}{\sum_{i=1}^n \vert p_{ji}^{\sigma_k}(t) \vert}
    ~~~k= 1,\dots,5,
\end{equation*}
where $p_{ji}^{\sigma_k}(t)$ are the entries of $P_{\sigma_k}(t)$ and 
$\sigma_k = \begin{bmatrix} 10 & 0.5 & 10 \end{bmatrix}^\top
+
(k-1) \begin{bmatrix} 2.5 & 0.25 & 7.5 \end{bmatrix}^\top$. 
In particular we have $\sigma_1 = \bar{\sigma}$.
Similarly, in \Cref{subf: genPrecA} we plot the associated generalized precisions 
\begin{equation*}
    p_k(t) = \det (P_{\sigma_k}(t)),~~~k = 1,\dots,5,
\end{equation*}
and in contrast to the first example we find no monotone dependence of the precision on the parameters.
We conclude that in this example the interpretation of the results is less straight forward and the system displays a structure noticeably different from the behavior of the scalar case. The energy minimizer $\widehat{x}_\mathrm{E}$ is particularly interesting as its dependence on the precision matrices becomes apparent. Finally, we point out that, as predicted in \Cref{cor: energyMin}, $\widehat{x}_\mathrm{E}$ offers the best approximation in terms of expected squared Mahalanobis distance, cf. \Cref{subf: expMah}.
%

%

%
\section{Conclusion}
%
In this work Kalman filter based concepts for the state estimation of uncertain linear systems are introduced, investigated, and numerically implemented. A particular focus lies on the minimizer of the expected energy which is characterized as the mean of the individual Kalman filters weighted with their precision matrices. It is theoretically proven and numerically illustrated 
that it minimizes the expected squared Mahalanobis distance to the Kalman filter arising from the realized, unknown parameter.
We find that in particular settings the described estimator is slightly outperformed by less elaborate approaches which is not surprising since it is defined as the minimizer of an averaged quantity. This issue might be addressed by considering more risk averse concepts optimizing metrics other than the average of the parameter family. 
Another interesting strain of research is the case of continuous parameter distributions such as normally distributed uncertain parameters. One could also investigate the impact of uncertainties on nonlinear systems leading to (higher order) extended Kalman filters or the so called Mortensen observer. 

\section*{Acknowledgement}
We thank P.~Guth (RICAM Linz) for many fruitful discussions on optimal control under uncertainty. 

%
%

\bibliographystyle{siam}
\bibliography{references} 

\begin{thebibliography}{10}

\bibitem{Ait36}
{\sc A.~C. Aitken}, {\em On least squares and linear combination of
  observations}, Proceedings of the Royal Society of Edinburgh, 55 (1936),
  p.~42–48.

\bibitem{BreKu21}
{\sc T.~Breiten and K.~Kunisch}, {\em Neural network based nonlinear
  observers}, System \& Control Letters, 148 (2021).

\bibitem{Cia78}
{\sc P.~G. Ciarlet}, {\em The Finite Element Method for Elliptic Problems},
  Studies in Mathematics and Its Applications, North-Holland, Amsterdam New
  York Oxford, first~ed., 1978.

\bibitem{DeGr70}
{\sc M.~H. DeGroot}, {\em Optimal statistical decisions}, McGraw-Hill, New York
  a.~o., first~ed., 1970.

\bibitem{DieEi94}
{\sc L.~Dieci and T.~Eirola}, {\em Positive definiteness in the numerical
  solution of {R}iccati differential equations}, Numerische Mathematik, 67
  (1994), p.~303–313.

\bibitem{Eve94}
{\sc G.~Evensen}, {\em Sequential data assimilation with a nonlinear
  quasi-geostrophic model using {M}onte {C}arlo methods to forecast error
  statistics}, J. Geophys. Res., 99 (1994), pp.~10143--10162.

\bibitem{Gir96}
{\sc V.~L. Girko}, {\em Random matrices}, in Handbook of Algebra,
  M.~Hazewinkel, ed., vol.~1 of Handbook of Algebra, North-Holland, 1996,
  pp.~27--78.

\bibitem{GutKuRo24}
{\sc P.~A. Guth, K.~Kunisch, and S.~S. Rodrigues}, {\em Tracking optimal
  feedback control under uncertain parameters}, Physica D: Nonlinear Phenomena,
  467 (2024), p.~134245.

\bibitem{Kal60}
{\sc R.~E. Kalman}, {\em A new approach to linear filtering and prediction
  problems}, Transactions of the ASME–Journal of Basic Engineering, 82
  (1960), pp.~35--45.

\bibitem{KalB61}
{\sc R.~E. Kalman and R.~S. Bucy}, {\em New results in linear filtering and
  prediction theory}, Transactions of the ASME–Journal of Basic Engineering,
  83 (1961), pp.~95--108.

\bibitem{KemScBy24}
{\sc H.~Kempka, C.~Schneider, and J.~Vybiral}, {\em Path regularity of the
  {B}rownian motion and the {B}rownian sheet}, Constructive Approximation, 59
  (2024), pp.~485--539.

\bibitem{KwaSiv72}
{\sc H.~Kwakernaak and R.~Sivan}, {\em Linear optimal control systems},
  Wiley-Interscience, a division of John Wiley {\&} Sons, Inc., New York,
  Chichester, Brisbane, Toronto, Singapore, first~ed., 1972.

\bibitem{Mah36}
{\sc P.~C. Mahalanobis}, {\em Reprint of: Mahalanobis, {P}.~{C}. (1936) "{O}n
  the generalised distance in statistics"}, The Indian Journal of Statistics,
  80 (2019), pp.~S1--S7.

\bibitem{Meh70}
{\sc R.~Mehra}, {\em On the identification of variances and adaptive {K}alman
  filtering}, IEEE Transactions on Automatic Control, 15 (1970), pp.~175--184.

\bibitem{Mor68}
{\sc R.~E. Mortensen}, {\em Maximum-likelihood recursive nonlinear filtering},
  Journal of Optimization Theory and Applications, 2 (1968), pp.~386--394.

\bibitem{Oks98}
{\sc B.~{\O}ksendal}, {\em Stochastic differential equations}, Universitext,
  Springer-Verlag, Berlin, fifth~ed., 1998.

\bibitem{PetMcf94}
{\sc I.~R. Petersen and D.~C. McFarlane}, {\em Optimal guaranteed cost control
  and filtering for uncertain linear systems}, IEEE Transactions on Automatic
  Control, 39 (1994), pp.~1971--1977.

\bibitem{PetSav99}
{\sc I.~R. Petersen and A.~V. Savkin}, {\em Robust Kalman Filtering for Signals
  and Systems with Large Uncertainties}, Birkh{\"a}user, Boston, first~ed.,
  1999.

\bibitem{Sch25}
{\sc J.~Schröder}, {\em Code for the paper "{D}eterministic {K}alman filters
  for uncertain dynamical systems"}, 2025.
\newblock doi:~10.5281/zenodo.15519854.

\bibitem{ShiJohMu07}
{\sc L.~Shi, K.~H. Johansson, and R.~M. Murray}, {\em Kalman filtering with
  uncertain process and measurement noise covariances with application to state
  estimation in sensor networks}, in 2007 IEEE International Conference on
  Control Applications, 2007, pp.~1031--1036.

\bibitem{Son98}
{\sc E.~D. Sontag}, {\em Mathematical control theory:~Deterministic finite
  dimensional systems}, Texts in Applied Mathematics 6, Springer, New York,
  second~ed., 1998.

\bibitem{Wie49}
{\sc N.~Wiener}, {\em Extrapolation, Interpolation, and Smoothing of Stationary
  Time Series}, The MIT Press, Cambridge, first~ed., 1949.

\bibitem{Wil04}
{\sc J.~C. Willems}, {\em Deterministic least squares filtering}, Journal of
  Econometrics, 118 (2004), pp.~341--373.

\bibitem{Xio08}
{\sc J.~Xiong}, {\em An introduction to stochastic filtering theory}, Oxford
  University Press, Oxford, first~ed., 2008.

\end{thebibliography}

\end{document}